\documentclass[12pt]{article}

\usepackage[english]{babel}
\usepackage{upref,amsfonts,amsxtra,a4wide,latexsym}
\usepackage{amssymb,a4wide,doublespace,amsmath,epsf,mathrsfs,theorem}
\usepackage[all]{xy}
 
\newtheorem{theorem}{Theorem}[subsection]
\newtheorem{lemma}[theorem]{Lemma}
\newtheorem{corollary}[theorem]{Corollary}

\newtheorem{proposition}[theorem]{Proposition}

{\theorembodyfont{\upshape}
\newtheorem{definition}[theorem]{Definition}
\newtheorem{remark}[theorem]{Remark}

\newtheorem{question}[theorem]{Question}}

\numberwithin{equation}{section}
\numberwithin{theorem}{section}

\newcommand{\image}{\mathrm{Im}}

\newcommand{\diag}{\mathrm{diag}}
\newcommand{\Id}{\mathrm{id}}

\newcommand{\Ad}{\mathrm{Ad}\,}

\newcommand{\ch}{{\rm {ch}}}

\newcommand{\e}{{\rm e}}

\newcommand{\cM}{{\cal M}}

\newcommand{\cN}{{\cal N}}

\newcommand{\cE}{{\cal E}}
\newcommand{\cO}{{\cal O}}

\newcommand{\C}{{\mathbb C}}
\newcommand{\Z}{{\mathbb Z}}
\newcommand{\N}{{\mathbb N}}

\newcommand{\T}{{\mathbb T}}

\newcommand{\F}{{\mathbb F}}
\newcommand{\cK}{{\cal K}}

\newcommand{\Cs}{{$C^*$-al\-ge\-bra}}

\newcommand{\gange}{\! \cdot \!}

\newcommand{\sh}{{$^*$-ho\-mo\-mor\-phism}}

\newcommand\eqdef{{\;\;\overset{\mbox{\scriptsize def}}{=}\;\;}}
\newcommand{\nprecsim}{{\operatorname{\hskip3pt \precsim\hskip-9pt |\hskip6pt}}}

\newenvironment{proof}[1][Proof:]%
{\begin{trivlist}\item[]\textbf{#1} }%
{\hbox{}\nobreak\hfill\quad\hbox{$\square$}\end{trivlist}}

\begin{document}
 
\title{A simple \Cs{} with a finite and an infinite projection}
\author{Mikael R\o rdam}
\date{}
\maketitle

\centerline{{\it{Dedicated to Richard V.\ Kadison on the occasion of
      his 75th birthday}}} 


\begin{abstract} \noindent An example is given of a simple, unital
  \Cs{} which contains an infinite and a non-zero finite projection. This
  \Cs{} is also an example of an infinite simple
  \Cs{} which is not purely infinite. A
  corner of this \Cs{} is a finite, simple, unital \Cs{} which is not
  stably finite. 

Our example shows that the type decomposition for von Neumann
factors does not carry over to simple \Cs s.

We also give an example of 
a simple, separable, \emph{nuclear}, \Cs{} in the UCT class which
contains an infinite and a non-zero finite projection. This nuclear
\Cs{} arises as a crossed product $D \rtimes_\alpha \Z$, where $D$ is
an inductive limit of type I \Cs s.
\end{abstract}

\section{Introduction}
\noindent The first interesting class of simple \Cs s (not counting
the simple von Neumann algebras) were the
UHF-algebras, also called Glimm algebras, constructed
by Glimm in 1959 (\cite{Gli:UHF}). Several other classes of simple \Cs
s were found over the following 25 years including the (simple)
AF-algebras,  the irrational rotation \Cs s, the free group \Cs s
$C^*_{\mathrm{red}}(\F_n)$ (and other reduced group \Cs s), the Cuntz
algebras $\cO_n$ and the Cuntz--Krieger algebras $\cO_A$, \Cs s
arising from minimal dynamical systems and from foliations, and
certain inductive limit \Cs s, among many other examples. Parallel with the
appearance of these examples of simple \Cs s it was asked if
there is a classification for simple \Cs s similar to the classification of
von Neumann factors into types. Inspired by work of Dixmier in the
1960's, Cuntz studied this and related questions about the structure
of simple \Cs s in his papers \cite{Cuntz:simple}, \cite{Cuntz:simpleII}, and 
\cite{Cuntz:dimension}. 

A von Neumann algebra is simple precisely when it is either a factor of
type I$_n$ for $n < \infty$ (in which case it is isomorphic to
$M_n(\C)$), a factor of type II$_1$, or a separable factor of type
III. This leads to the question if (non type I) simple \Cs s can
be divided into two subclasses, one that resembles type II$_1$ factors 
and another that resembles type III factors. A II$_1$ 
factor is an infinite dimensional factor in which all projections are
finite (in the sense of Murray--von Neumann's comparison theory for
projections), and II$_1$ factors have a unique trace. A factor is of
type III if all its non-zero projections are infinite, and type III
factors admit no traces. Cuntz asked in \cite{Cuntz:simpleII} if each
simple \Cs{} similarly must have the property that its (non-zero)
projections either all are finite or all are infinite. Or can a simple
\Cs{} contain both a (non-zero) finite and an infinite projection? We
answer the latter question in the affirmative. In other words, we
exhibit a simple (non type I) \Cs s that neither corresponds to a type
II$_1$ or to a type III factor.  

It was shown in the early 1980's that simple \Cs s, in contrast to von
Neumann factors, can fail to have non-trivial projections. Blackadar
(\cite{Bla:projless}) and Connes (\cite{Con:Thom}) found examples of
unital, simple \Cs s  with no
projections other than $0$ and $1$---before it was shown that
$C^*_{\mathrm{red}}(\F_2)$ is a simple unital \Cs{} with no
non-trivial projections. Simple \Cs s can fail to have projections in
a more severe way: Blackadar found in \cite{Bla:projless2} an example
of a stably projectionless simple \Cs s. (A 
\Cs{} $A$ is stably projectionless if $0$ is the only projection in $A
\otimes \cK$.) Blackadar and Cuntz proved in \cite{BlaCuntz:infproj}
that every stably projectionless simple \Cs{} is finite in the sense
of admitting a (densely defined) quasitrace. (Every quasitrace on an
exact \Cs{} 
extends to a trace as shown by Haagerup \cite{Haa:quasi} (and
Kirchberg \cite{Kir:quasitraces}).) These results lead to the dichotomy
for a
simple \Cs{} $A$: Either $A$ admits a (densely defined) quasitrace
(in which case $A$ is stably finite), or $A$ is \emph{stably
  infinite}, i.e., $A \otimes \cK$ contains an infinite
projection. 

Cuntz defined in \cite{Cuntz:KOn} a simple \Cs{} to be \emph{purely
  infinite} if all its non-zero here\-di\-tary sub-\Cs s contain an 
infinite projection. Cuntz showed in \cite{Cuntz:On} that his
algebras $\cO_n$, $2 \le n \le \infty$, are simple and purely
infinite. The separable, nuclear, simple, purely infinite \Cs  s are
classified up to isomorphism by $K$- or $KK$-theory by the spectacular 
theorem of Kirchberg (\cite{Kir:Michael}
and \cite{Kir:fields}) and Phillips (\cite{Phi:class}). This result
has made it an important
question to decide which simple \Cs s are purely infinite. We show here
that not all stably infinite simple \Cs s $A$ are purely
infinite. 

Villadsen (\cite{Vil:perforation}) was the first to show that the
$K_0$-group of a simple \Cs{} need not be weakly
unperforated; Villadsen (\cite{Vil:sr=n}) also showed that a unital,
finite, simple \Cs{} can have stable rank different from one---thus
answering in the negative two longstanding open questions for simple
\Cs s.

If $B$ is a unital, simple \Cs{} with an infinite and a non-zero
finite projection, then its semigroup
of Murray--von Neumann equivalence classes of 
projections must fail to be weakly unperforated (see Remark~\ref{rem:comp}).
It is therefore no surprise that Villadsen's ideas play a crucial role
in this article. Our article is also a continuation of the work by the  
author in \cite{Ror:sns} and \cite{Ror:sums} where it is shown that
one can find a \Cs{} $A$ such that $M_2(A)$ is stable but $A$ is
not stable; and, related to this, one can find a (non-simple)
unital \Cs{} $B$, such that $B$ is finite and $M_2(B)$ is
properly infinite. We show here (Theorem~\ref{thm:B})
that one can make this example simple
by passing to a suitable inductive limit. 

In Section~6 (added March 2002) an example is given of a crossed
product \Cs{} $D \rtimes_\alpha \Z$, where $D$ is an inductive limit
of type I \Cs s, such that $D \rtimes_\alpha \Z$ is simple and
contains an infinite and a non-zero finite projection. This new
example is nuclear and separable. It shows that simple \Cs s with this
rather pathological behavior can arise from a quite natural
setting. It shows that Elliott's classification conjecture (in its present
formulation) does not hold (cf.\ Corollary~\ref{cor:conjecture}); and it
also serves as an example of a separable nuclear simple \Cs{} that is
tensorially prime (cf.\ Corollary~\ref{cor:prime}).

I thank Bruce Blackadar, Joachim Cuntz, George Elliott, and Eberhard
Kirchberg for valuable discussions and for their comments to
earlier versions of this 
manuscript. I thank Paul M.\ Cohn and Ken Goodearl for explaining the example
included in Remark~\ref{rm:ring}.  I also thank the referee for
suggesting several improvements to this article (including a significant
simplification of  
Proposition~\ref{prop:varphi}~(ii) and (iii)). 

This work was done in the spring of 2001 while the author visited the
University of California, Santa Barbara. I thank Dietmar Bisch for
inviting me and for his warm hospitality. 

The present revised version (with the nuclear example in Section~6 and
where the construction in Section~5 is simplified) was completed in
March 2002. A part of the work leading to this construction was
obtained during a visit in
January 2002 to the University of M\"unster. I thank Joachim Cuntz and
Eberhard Kirchberg for their hospitality, and I am indebted to
Eberhard Kirchberg for several conversations during the visit that
led me to this construction.

\section{Finite, infinite, and properly infinite projections}
\label{sec:prelim} 

\noindent A projection $p$ in a \Cs{} $A$ is called \emph{infinite} if it is
equivalent (in the sense of Murray and von Neumann) to a proper
subprojection of itself; and $p$ is said to be \emph{finite}
otherwise. If $p$ is non-zero and if there are mutually orthogonal 
subprojections $p_1$ and $p_2$ of $p$ such that $p \sim p_1 \sim p_2$, 
then $p$ is \emph{properly infinite}. A unital \Cs{} is said to be
\emph{properly infinite} if its unit is a properly infinite projection.

If $p$ and $q$ are projections in $A$, then let $p \oplus q$ denote
the projection $\diag(p,q)$ in $M_2(A)$. Two
projections $p \in M_n(A)$ and $q \in M_m(A)$ can be compared as
follows: Write $p \sim q$ if there exists $v$ in $M_{m,n}(A)$ such
that $v^*v = p$ and $vv^* = q$, and write $p \precsim q$ if $p$ is
equivalent (in this sense) to a subprojection of $q$. 

In the proposition below, where some well-known properties of
properly infinite projections are recorded, $\cO_\infty$ denotes the Cuntz algebra
generated by infinitely many isometries with pairwise orthogonal range 
projections, and $\cE_2$ is the Cuntz--Toeplitz algebra generated by two
isometries with orthogonal range projections (\cite{Cuntz:On}). 

\begin{proposition} \label{prop:prop}
The following five conditions are equivalent for every non-zero
projection $p$ in a \Cs{} $A$:
\begin{enumerate}
\item $p$ is properly infinite;
\item $p \oplus p \precsim p$;
\item there is a unital \sh{} $\cE_2 \to pAp$;
\item there is a unital \sh{} $\cO_\infty \to pAp$;
\item for every closed two-sided ideal $I$ in $A$, either $p \in I$ or
  $p+I$ is infinite in $A/I$.
\end{enumerate}
\end{proposition}

\noindent The equivalences between (i), (ii), and (iii) are
trivial. The equivalence between (iii) and (iv) follows from the fact
that there are unital embeddings $\cE_2 \to \cO_\infty$ and
$\cO_\infty \to \cE_2$. The equivalence between (i) and (v) is proved in
  \cite[Corollary 3.15]{KirRor:pi}; a result that extends Cuntz'
  important observation from \cite{Cuntz:simple} that every infinite projection
  in a simple \Cs{} is properly infinite.

We shall use the following two well-known results about properly
infinite projections. 

\begin{lemma} \label{lm:prop}
Let $p$ and $q$ be projections in a \Cs{} $A$. Suppose that $p$ is
properly infinite. Then $q \precsim p$ if and only if $q$ belongs to
the closed two-sided ideal in $A$ generated by $p$.
\end{lemma}

\begin{proof} If $q \precsim p$, then, by definition, $q \sim q_0 \le
  p$ for some projection $q_0$ in $A$. This entails that $q$
  belongs to the ideal generated 
  by $p$. Conversely, if $q$ belongs to the ideal generated by $p$,
  then $q \precsim \bigoplus_{j=1}^n p$ for some $n$ (cf.\
  \cite[Exercise 4.8]{RorLarLau:k-theory}), and $\bigoplus_{j=1}^n p
  \precsim p$ if $p$ is properly infinite by iterated applications of
  Proposition~\ref{prop:prop}~(ii). 
\end{proof}

\begin{proposition} \label{prop:inductive-limit}
Let $B$ be the inductive limit of a sequence $B_1 \to B_2 \to B_3 \to
\cdots$ of unital \Cs s with unital connecting maps. Then $B$
is properly infinite if and only if $B_n$ is properly infinite for all $n$
larger than some $n_0$.
\end{proposition}

\begin{proof} If $B_n$ is properly infinite for some $n$, then there
  are unital \sh s $\cE_2 \to B_n \to B$, and hence $B$ is
  properly infinite. Conversely, if $B$ is properly infinite, then
  there is a
  unital \sh{} $\cE_2 \to B$. The \Cs{} $\cE_2$ is semiprojective, as
  shown by Blackadar in \cite{Bla:shape}. By semiprojectivity (see
  again \cite{Bla:shape}), the
  unital \sh{} $\cE_2 \to B$ lifts to a unital \sh{} $\cE_2
  \to \prod_{n=n_0}^\infty B_n$ for some $n_0$. This shows that
  $B_n$ is properly infinite for all $n \ge n_0$.
\end{proof}

\section{Vector bundles over products of spheres} \label{sec:bundle}

\noindent We consider here complex vector bundles over the
sphere $S^2$ and over finite products of spheres, $(S^2)^n$. 

For each $k \le n$, let $\pi_k \colon (S^2)^n \to S^2$ denote the $k$th 
coordinate mapping, and let $\rho_{m,n} \colon (S^2)^{m} \to (S^2)^n$ be
given by
\begin{equation} \label{eq:rho_{m,n}}
\rho_{m,n}(x_1,x_2, \dots, x_{m}) = (x_1,x_2, \dots, x_n), \qquad
(x_1,x_2, \dots, x_m) \in (S^2)^m.
\end{equation}
when $m \ge n$.

Whenever $f \colon X \to Y$ is a continuous map and $\xi$ is a
$k$-dimensional complex vector bundle over $Y$, let $f^*(\xi)$ denote
the vector bundle over $X$ induced by $f$. Let $e(\xi) \in H^{2k}(Y,\Z)$
denote the Euler class of $\xi$. Denote also by $f^*$ the induced map
$H^*(Y,\Z) \to H^*(X,\Z)$. By functoriality of the Euler class we have 
$f^*(e(\xi)) = e(f^*(\xi))$. 

For any vector bundle $\xi$ over $(S^2)^n$ and for every $m \ge n$ we
have a vector bundle $\xi'=\rho_{m,n}^*(\xi)$ over $(S^2)^m$. It
follows from the K\"unneth Theorem (see \cite[Theorem
A6]{MilSta:classes}), that the map  
$$\rho_{m,n}^* \colon H^*((S^2)^n,\Z) \to H^*((S^2)^m,\Z)$$
is injective; so if $e(\xi)$ is non-zero, then so is $e(\xi')$. Our
main concern with vector bundles will be whether or not they have
non-zero Euler class, and from that point of view it does not matter
if we replace the base space $(S^2)^n$ with $(S^2)^m$ for
some $m \ge n$. 

We remind the reader of some properties of the Euler class for complex
vector
bundles $\xi_1,\xi_2, \dots, \xi_n$ over a base space $X$. First of
all we have the product formula (see \cite[Property 9.6]{MilSta:classes}):
\begin{equation} \label{eq:Euler1}
e(\xi_1 \oplus \xi_2 \oplus \cdots \oplus \xi_n) = e(\xi_1) \gange
e(\xi_2) \cdots  e(\xi_n).
\end{equation}
Let $\theta$  denote the trivial complex line bundle over 
$X$. The Euler class of
$\theta$ is zero; and so it follows from the product formula that
$e(\xi)=0$ whenever $\xi$ is a complex vector 
bundle that dominates $\theta$ in the sense that
$\xi \cong \theta \oplus \eta$ for some complex vector bundle $\eta$. 

Combining the formula
$$\ch(\xi) = 1 + e(\xi) + \frac{1}{2}e(\xi)^2 + \frac{1}{6}\e(\xi)^3 + 
\cdots,$$
that relates the Chern character and the Euler class of a complex line
bundle $\xi$ (see \cite[Problem 16-B]{MilSta:classes}), with the fact
that the Chern character is multiplicative, yields the formula
\begin{equation} \label{eq:Euler2}
e(\xi_1 \otimes \xi_2 \otimes \cdots \otimes \xi_n) = e(\xi_1) +
e(\xi_2) + \cdots + e(\xi_n),
\end{equation}
that holds for all complex \emph{line bundles} $\xi_1, \dots, \xi_n$ over
$X$.

{\bf{Let $\zeta$ be a complex line bundle over $S^2$ such that
its Euler class $e(\zeta)$, which is an element in $H^2(S^2,\Z)$, is
non-zero.}} (Any such line bundle will do, 
but the reader may take $\zeta$ to be the Hopf bundle over $S^2$.) 
For each natural number $n$ and for each non-empty, finite subset $I =
\{n_1,n_2, \dots, n_k\}$ of $\N$ define complex line bundles
$\zeta_n$ and $\zeta_I$ over $(S^2)^m$ (for all $m \ge n$,
respectively, $m \ge \max\{n_1,\dots,n_k\}$) by
\begin{equation} \label{eq:zeta_n}
\zeta_n = \pi_n^*(\zeta), \qquad  \zeta_I = \zeta_{n_1} \otimes
\zeta_{n_2} \otimes \cdots \otimes \zeta_{n_k},
\end{equation}
where, as above, $\pi_n \colon (S^2)^m \to S^2$ is the $n$th coordinate
map. The Euler classes (in $H^2((S^2)^m,\Z)$) of these line bundles 
are by functoriality and equation~\eqref{eq:Euler2} given by 
\begin{eqnarray} \label{eq:Euler3}
e(\zeta_n) & = &\pi_n^*(e(\zeta)), \\ \label{eq:Euler4}
e(\zeta_I) & = & \pi_{n_1}^*(e(\zeta))+\pi_{n_2}^*(e(\zeta))+ \cdots
+\pi_{n_k}^*(e(\zeta)).
\end{eqnarray}

\begin{lemma} \label{lm:2zeta}
For each $n$ and for each $m \ge n$ there is a complex line bundle
$\eta_n$ over $(S^2)^m$ such that $\zeta_n \oplus \zeta_n \cong \theta
\oplus \eta_n$. 
\end{lemma}

\begin{proof} Since 
$$\dim(\zeta \oplus \zeta) = 2 > 1 \ge \textstyle{\frac{1}{2}}(\dim(S^2)-1),$$
it follows from \cite[9.1.2]{Hus:fibre} that there 
  is a complex vector bundle $\eta$ over $S^2$ of dimension
  $\dim(\eta) = 2-1=1$ such that $\zeta \oplus \zeta \cong \theta
  \oplus \eta$. We conclude that
$$\zeta_n \oplus \zeta_n = \pi_n^*(\zeta \oplus \zeta) \cong
\pi_n^*(\theta \oplus \eta) = \theta \oplus \pi_n^*(\eta).$$
\end{proof}

\begin{proposition} \label{prop:marriage1}
Let $I_1, I_2, \dots, I_m$ be non-empty, finite subsets of $\N$. The
following three conditions are equivalent:
\begin{enumerate}
\item $e(\zeta_{I_1} \oplus \zeta_{I_2} \oplus \cdots \oplus
  \zeta_{I_m}) \ne 0$.
\item For all subsets $F$ of $\{1,2, \dots, m\}$ we have $\big|
  \bigcup_{j \in F} I_j \big| \ge |F|$.
\item There exists a matching $t_1 \in I_1, t_2 \in
  I_2, \dots, t_m \in I_m$ (i.e., the elements $t_1, \dots, t_m$ are
  pairwise distinct).
\end{enumerate}
\end{proposition}

\begin{proof} Choose $N$ large enough so that each $\zeta_{I_j}$ is
  a vector bundle over $(S^2)^N$.

(ii) $\Leftrightarrow$ (iii) is the Marriage Theorem (see any textbook on
combinatorics). 

(i) $\Rightarrow$ (ii). Assume that  $\big| \bigcup_{j \in F} I_j
\big| < |F|$ for some (necessarily non-empty) subset $F = \{j_1, j_2, \dots, j_k\}$ of
$\{1,2, \dots, m\}$, and write
$$J \eqdef \bigcup_{j \in F} I_j =
\{n_1,n_2, \dots, n_l\}.$$ 
Let $\rho \colon (S^2)^N \to (S^2)^l$ be given by
$\rho(x) = (\pi_{n_1}(x), \pi_{n_2}(x), \dots, \pi_{n_l}(x))$.
Then 
$$\xi \eqdef \zeta_{I_{j_1}} \oplus \zeta_{I_{j_2}} \oplus \cdots
\oplus \zeta_{I_{j_k}} = \rho^*(\eta)$$
for some $k$-dimensional vector bundle $\eta$ over
$(S^2)^l$. Now, $e(\eta)$ belongs to $H^{2k}((S^2)^l,\Z)$, and
$H^{2k}((S^2)^l,\Z)=0$  because $2k > 2l$. Hence $e(\xi) =
\rho^*(e(\eta)) =0$, so by the  
product formula \eqref{eq:Euler1} we get
$$e(\zeta_{I_1} \oplus \zeta_{I_2} \oplus \cdots \oplus \zeta_{I_m}) = 
e(\xi) \gange \prod_{j \notin F} e(\zeta_{I_j}) = 0.$$

(iii) $\Rightarrow$ (i). Put 
$$x_j = \pi_j^*(e(\zeta)) \in H^2((S^2)^N,\Z), \qquad j=1,2, \dots, N.$$ 
The element
$$z = x_1 \gange x_2  \cdots  x_N \in H^{2N}((S^2)^N,\Z)$$
is non-zero by the K\"unneth Theorem (\cite[Theorem
A6]{MilSta:classes}). Using that $x_i^2 =0$ and that
$x_ix_j=x_jx_i$ for all $i,j$ it follows that if $i_1,i_2, \dots, i_N$ 
belong to $\{1,2,\dots,N\}$, then
\begin{equation} \label{eq:product}
x_{i_1} \gange x_{i_2}  \cdots x_{i_N} = \begin{cases} z, 
  & \text{if} \; \; i_1, \dots, i_N \; \; \text{are distinct,}\\ 0, &
  \text{otherwise.} \end{cases}
\end{equation}
Now, by \eqref{eq:Euler1} and \eqref{eq:Euler4}, 
\begin{eqnarray*}
e(\zeta_{I_1} \oplus \zeta_{I_2} \oplus \cdots \oplus \zeta_{I_m}) &=&
e(\zeta_{I_1}) \gange  e(\zeta_{I_2}) \cdots e(\zeta_{I_m})  
\\ & = & \big(\sum_{i \in I_1} \pi_i^*(e(\zeta))\big) \gange
\big(\sum_{i \in I_2} \pi_i^*(e(\zeta))\big) \cdots \big(\sum_{i \in I_m}
\pi_i^*(e(\zeta))\big)  \\
&=& \big(\sum_{i \in I_1} x_i \big) \gange
\big(\sum_{i \in I_2} x_i \big) \cdots \big(\sum_{i \in I_m}
x_i \big)  \\
&=& \sum_{(i_1, \dots, i_m) \in I_1 \times \dots \times I_m} x_{i_1} \gange
  x_{i_2} \cdots x_{i_m}.
\end{eqnarray*}
Assume that (iii) holds, and write
$$\{1,2, \dots, N\} \setminus \{t_1,t_2, \dots, t_m\} = \{s_1, s_2,
\dots, s_{N-m}\}.$$
Let $k$ denote the number of permutations $\sigma$ on $\{1,2, \dots,
m\}$ such that $t_{\sigma(j)} \in I_j$ for $j=1,2, \dots, m$. The
identity permutation has this property, so $k \ge 1$. The formula for
$e(\zeta_{I_1} \oplus \cdots  
\oplus \zeta_{I_m})$ above and equation \eqref{eq:product} yield
$$e(\zeta_{I_1} \oplus \zeta_{I_2} \oplus \cdots \oplus
\zeta_{I_m})\gange x_{s_1} \gange 
x_{s_2} \cdots x_{s_{N-m}} = kz \ne 0.$$ 
It follows that $e(\zeta_{I_1} \oplus \cdots \oplus \zeta_{I_m}) \ne 0$ as desired.
\end{proof}

\section{Projections in a certain multiplier algebra} \label{sec:multiplier}

\noindent There is a well-known one-to-one correspondence between 
isomorphism classes of complex vector bundles over a compact Hausdorff
space $X$ 
and Murray--von Neumann equivalence classes of projections in matrix
algebras over $C(X)$ (and in $C(X) \otimes \cK$). The vector bundle
corresponding to a projection $p$ in $M_n(C(X)) =  C(X, M_n(\C))$ is 
$$\xi_p = \{(x,v) : x \in X, \; v \in p(x)(\C^n) \},$$
(equipped with the topology given from the natural inclusion $\xi_p
\subseteq X \times \C^n$),
so that the fibre $(\xi_p)_x$ over $x \in X$ is the range of the projection
$p(x)$. If $p$ and $q$ are two projections in $C(X) \otimes \cK$, then 
$\xi_p \cong \xi_q$ if and only if $p \sim q$. It follows from Swan's
theorem, which to each complex vector bundle $\xi$ gives a complex
vector bundle  
$\eta$ such that $\xi \oplus \eta$ is isomorphic to 
the trivial $n$-dimensional complex vector bundle over $X$ for some
$n$, that every complex vector bundle is isomorphic to $\xi_p$ for
some projection $p$ in $M_n(C(X))$ for some $n$. 

View each matrix
algebra $M_n(\C)$ as a sub-\Cs{} of $\cK$ via the embeddings
$$\xymatrix{\C \, \ar@{^(->}[r]& M_2(\C) \,  \ar@{^(->}[r] & M_3(\C) \,
  \ar@{^(->}[r] & \cdots \ar@{^(->}[r]  & \cK,}$$
where $M_n(\C)$ is mapped into the upper left corner of
$M_{n+1}(\C)$. Identify $C(X,\cK)$ with $C(X)
\otimes \cK$ and identify $C(X,M_n(\C))$ with $C(X) \otimes M_n(\C)$.

In Section~\ref{sec:bundle} we picked a non-trivial complex line
bundle $\zeta$ over $S^2$ (eg., the Hopf bundle). This line bundle $\zeta$
corresponds to a projection $p$ in some matrix algebra over $C(S^2)$,
and, as is well known, such a projection $p$ can be found in
$M_2(C(S^2)) = C(S^2,M_2)$. (The projection $p \in M_2(S^2,M_2)$ corresponding to the
Hopf bundle is in operator algebra texts often referred to as the Bott
projection.)
Put
$$Z= \prod_{n=1}^\infty S^2.$$ 
Let $\pi_n \colon Z \to S^2$ be the $n$th coordinate map, and let
$\rho_{\infty,n} \colon Z \to (S^2)^n$ be given by
$$\rho_{\infty,n}(x_1,x_2,x_3, \dots) = (x_1,x_2, \dots, x_n), \qquad
(x_1,x_2, x_3, \dots) \in Z.$$
With ${\widehat{\rho}_n} \colon
C((S^2)^n) \to C((S^2)^{n+1})$ being the \sh{} induced by the map
$\rho_n = \rho_{n+1,n}$ defined in \eqref{eq:rho_{m,n}} we obtain that 
$C(Z)$ is the inductive limit
$$\xymatrix{C(S^2)  \ar[r]^-{\widehat{\rho}_1} & C((S^2)^2) 
  \ar[r]^-{\widehat{\rho}_2} & C((S^2)^3)  \ar[r]^-{\widehat{\rho}_3} &
  \cdots \ar[r] & C(Z)}$$
with inductive limit maps ${\widehat{\rho}}_{\infty,n} \colon
  C((S^2)^n) \to C(Z)$.

For $n$ in $\N$ and for each
non-empty finite subset $I = \{n_1,n_2, \dots, n_k\}$ of $\N$, let
$p_n$ and $p_I$ be the projections in $C(Z) \otimes \cK = C(Z,\cK)$ given by
\begin{alignat}{1}
\label{eq:p_n} p_n(x) & \; = \; p(x_n), \\
\begin{split}
\label{eq:p_I} p_I(x) & \; = \; p(x_{n_1}) \otimes p(x_{n_2}) \otimes \cdots
\otimes p(x_{n_k}) \\
& \; = \; p_{n_1}(x) \otimes p_{n_2}(x) \otimes \cdots \otimes p_{n_k}(x),
\end{split}
\end{alignat}
for all $x =(x_1,x_2, \dots) \in Z$ (identifying $M_2$, respectively,
$M_2 \otimes M_2 \otimes \cdots \otimes M_2$ with sub-\Cs s of $\cK$).

We shall now make use of the multiplier algebra, $\cM(C(Z)
\otimes \cK)$, of $C(Z) \otimes \cK = C(Z,\cK)$. We can identify this
multiplier algebra with the set of all bounded functions $f \colon Z
\to B(H)$ for which $f$ and $f^*$ are continuous, when $B(H)$, the
bounded operators on the Hilbert space $H$ on which $\cK$ acts, is
given the strong operator topology.

It is convenient to have a convention for adding finitely or 
infinitely many projections in $\cM(C(Z)
\otimes \cK)$, or more generally in $\cM(A)$, where $A$ is any stable
\Cs{}---a convention that extends the notion of forming direct sums of
projections discussed in Section~\ref{sec:prelim}. 

Assuming that $A$ is a stable \Cs, so that $A = A_0 \otimes \cK$ for some \Cs{}
$A_0$, then we can
take a sequence $\{T_j\}_{j=1}^\infty$ of
isometries in 
$\C \otimes B(H) \subseteq \cM(A_0 \otimes \cK) = \cM(A)$ such that $ 1 =
\sum_{j=1}^\infty T_jT_j^*$ 
in the strict topology. (Notice that $1$ is a properly
infinite projection in $\cM(A)$.)
For any sequence $q_1,q_2, \dots$ of
projections in $A$ and for any sequence $Q_1,Q_2, \dots$ of
projections in $\cM(A)$, define
\begin{eqnarray}
\label{eq:DS1} q_1 \oplus q_2 \oplus \cdots \oplus q_n &=&  
\sum_{j=1}^n T_jq_jT_j^* \; \in \; A,\\
\label{eq:DS2} \bigoplus_{j=1}^\infty q_j & = & \sum_{j=1}^\infty
T_jq_jT_j^* \; \in \; \cM(A),\\
\label{eq:DS3} Q_1 \oplus Q_2 \oplus \cdots \oplus Q_n & = &
\sum_{j=1}^n T_jQ_jT_j^* \; \in \; \cM(A),\\
\label{eq:DS4} \bigoplus_{j=1}^\infty Q_j & = &
\sum_{j=1}^\infty T_jQ_jT_j^* \; \in \; \cM(A),
\end{eqnarray}
Observe that $q_j'=T_jq_jT_j^* \sim q_j$, that the projections
$q_1',q_2', \dots$ are mutually orthogonal, and that the sum
$\sum_{j=1}^\infty q_j'$ is strictly convergent. The projections in
\eqref{eq:DS1}--\eqref{eq:DS4} are, up to unitary equivalence in
$\cM(A)$, independent of the choice of isometries
$\{T_j\}_{j=1}^\infty$. Indeed, if $\{R_j\}_{j=1}^\infty$ is another
sequence of isometries in $\cM(A)$ with  $ 1 =
\sum_{j=1}^\infty R_jR_j^*$, then $U = \sum_{j=1}^\infty R_jT_j^*$ is
a unitary element in $\cM(A)$ and
$$\sum_{j=1}^\infty R_j X_j R_j^* = U \big( \sum_{j=1}^\infty
T_jX_jT_j^* \big) U^*$$
for any bounded sequence $\{X_j\}_{j=1}^\infty$ in $\cM(A)$.
It follows in particular that
\begin{equation} \label{eq:permute}
\bigoplus_{j=1}^\infty q_j \sim \bigoplus_{j=1}^\infty q_{\sigma(j)}
\end{equation}
for every permutation $\sigma$ on $\N$.

In the lemma below the correspondence between projections and vector
bundles is the mapping $p \mapsto \xi_p$ defined at the beginning of
this section. By identifying the projections $p_n, p_I,
p_{I_1}, \dots, p_{I_k}$ with projections in $C((S^2)^N) \otimes \cK$,
where $N$ is 
any integer large enough to ensure that these projections belong to the image
of $$\widehat{\rho}_{\infty,N} \otimes \Id_{\cK} \colon C((S^2)^N)
\otimes \cK \to C(Z) 
\otimes \cK,$$ 
we can take the base space to be $(S^2)^N$.

\begin{lemma} \label{lm:Euler-p}  Let 
$\zeta_n$ and $\zeta_I$ be the complex line bundles defined
in \eqref{eq:zeta_n}.
\begin{enumerate}
\item The vector bundle $\zeta_n$ corresponds to $p_n$ for each $n$ in $\N$.
\item The vector bundle $\zeta_I$ corresponds to $p_I$ for each
  non-empty finite subset $I$ of  $\N$.
\item The vector bundle  $\zeta_{I_1}  
  \oplus \zeta_{I_2} \oplus \cdots \oplus \zeta_{I_k}$ corresponds to
  $p_{I_1} \oplus p_{I_2} 
  \oplus  \cdots \oplus p_{I_k}$ whenever $I_1,  
  \dots, I_k$ are non-empty finite subsets of $\N$.
\end{enumerate}
\end{lemma}

\begin{proof} 
(i). Since $p$ corresponds to $\zeta$, $p_n = p \circ \pi_n$
  corresponds to $\zeta_n = \pi_n^*(\zeta)$, where $\pi_n \colon
  (S^2)^N \to S^2$ is the $n$th coordinate map. 

(ii). Write $I = \{n_1,n_2,\dots,n_k\}$. We shall here view $p_n$ as a 
projection in $C((S^2)^N,M_2)$ and $p_I$ as a projection in
$C((S^2)^N,M_2 \otimes \cdots \otimes M_2)$. 
By (i), $\zeta_n$ is the complex line bundle over $(S^2)^N$ whose
fibre over $x \in (S^2)^N$ is equal to $p_n(x)(\C^2)$.  The fibre of
the complex line bundle $\zeta_I = \zeta_{n_1} \otimes
\zeta_{n_2} \otimes \cdots \otimes \zeta_{n_k}$ over $x \in (S^2)^N$ is
by definition
\begin{eqnarray*}
(\zeta_I)_x & = & (\zeta_{n_1})_x \otimes (\zeta_{n_2})_x \otimes
\cdots \otimes (\zeta_{n_1})_x \\  & = & 
p_{n_1}(x)(\C^2) \otimes p_{n_2}(x)(\C^2) \otimes \cdots \otimes
p_{n_k}(x)(\C^2) \\ & = & p_I(x)(\C^2 \otimes \C^2 \otimes \cdots \otimes
\C^2).
\end{eqnarray*}
This shows that $\zeta_I$ corresponds to $p_I$.

(iii). This follows from (ii) and additivity of the map $p \mapsto \xi_p$.
\end{proof}

\noindent The two next results are formulated for an arbitrary stable
\Cs{} $A$ and its multiplier algebra $\cM(A)$, but it shall primarily
be used in the case where $A = C(Z) \otimes \cK$.

The lemma below is a trivial, but much used, generalization of
\eqref{eq:permute}: 

\begin{lemma} \label{lm:rearranging}
Let $A$ be a stable \Cs, and let $q_1,q_2, \dots$ and $r_1, r_2,
\dots$ be two sequences of 
projections in $A$. Assume that there is a permutation
$\sigma$ on $\N$ such that $q_j \precsim
r_{\sigma(j)}$, respectively $q_j \sim
r_{\sigma(j)}$, in $A$ for all $j$ in $\N$. Then
$\bigoplus_{j=1}^\infty q_j \precsim \bigoplus_{j=1}^\infty r_j$,
respectively $\bigoplus_{j=1}^\infty q_j \sim \bigoplus_{j=1}^\infty r_j$,
in $\cM(A)$.
\end{lemma}

\noindent We shall also use the following two well-known facts about
projections in multiplier algebras. An element in a \Cs{} $A$
is said to be \emph{full in $A$} if it is not contained in any proper closed
two-sided ideal of $A$.

\begin{lemma} \label{lm:Qsim1}
Let $A$ be a stable \Cs{}. The following three conditions are
equivalent for all projections $Q$ in $\cM(A)$:
\begin{center}
{\rm{(i)}} $Q \sim 1$, \hspace{.5cm} {\rm{(ii)}} $Q$ is
properly infinite and full in $\cM(A)$, \hspace{.5cm} {\rm{(iii)}}
$1 \precsim Q$. 
\end{center}
\end{lemma}

\begin{proof}
(i) $\Rightarrow$ (iii) is trivial. Assume that $1 \precsim Q$. Then
$Q$ is full in $\cM(A)$ (the closed two-sided ideal in $\cM(A)$
generated by $Q$ contains $1$ and hence all of $\cM(A)$). It was noted
above \eqref{eq:DS1} that $1$ is properly infinite in $\cM(A)$, and so
$Q \oplus Q \le 1 \oplus 1 \precsim 1 \precsim Q$, whence $Q$ is
properly infinite; cf.\ Proposition~\ref{prop:prop}. This proves (iii)
$\Rightarrow$ (ii). Assume finally that $Q$ is properly infinite and
full in $\cM(A)$. Since $K_0(\cM(A))=0$ (see \cite[Proposition
12.2.1]{Bla:k-theory}) the two projections $Q$ and $1$ represent the same
element in $K_0(\cM(A))$; and since these two projections both are
properly infinite and full they must be Murray--von Neumann equivalent
(see \cite[Section 1]{Cuntz:KOn} or \cite[Exercise 4.9
(iii)]{RorLarLau:k-theory}), i.e., $Q \sim 1$.  
\end{proof}

\begin{lemma} \label{lm:q<Q}
Let $A$ be a stable \Cs{} and let $q,q_1,q_2, \dots$ be projections in
$A$. If $q \precsim \bigoplus_{j=1}^\infty q_j$ in 
  $\cM(A)$, then $q \precsim q_1 \oplus q_2 \oplus
  \cdots \oplus q_k$ in $A$ for some $k$. 
\end{lemma}

\begin{proof} We have $\bigoplus_{j=1}^\infty q_j = \sum_{j=1}^\infty
  q_j' \; (= Q)$ for some strictly summable sequence of mutually
  orthogonal projections 
  $q_1',q_2', \dots$ in $A$ with $q_j' \sim q_j$. By
  the assumption that $q \precsim Q$ there is a partial isometry $v$
  in $\cM(A)$ such that $vv^* = q$ and $v^*v
  \le Q$. As $v = qv$, $v$ belongs to $A$, and by the strict
  convergence of the sum $Q = \sum_{j=1}^\infty q_j'$ there is $k$ such that
$$\|v - v\sum_{j=1}^k q_j'\| < 1/2.$$
Put $x = v\sum_{j=1}^k q_j'$. Then $xx^* \le q$, $x^*x
\le q_1' + \cdots + q_k'$, and $\|xx^*-q\| <1$. This shows that $xx^*$
is invertible in $qAq$ with inverse
$(xx^*)^{-1}$. Put $u = (xx^*)^{-1/2}x$. Then $uu^*=q$ and $u^*u \le
q_1' + \cdots + q_k'$, whence $q \precsim  q_1  \oplus
  \cdots \oplus q_k$.
\end{proof}

\noindent Let $g$ be a constant one-dimensional projection
  in $C(Z,\cK)=C(Z) \otimes \cK$ (that corresponds to the trivial complex line
  bundle $\theta$ over $X$). The (easy-to-prove) statement in part (iii) of the
  proposition below is not used in this paper, but it may have some
  independent interest.

\begin{proposition} \label{prop:marriage2} Let $I_1,I_2, \dots$ be a
  sequence of non-empty,
  finite subsets of $\N$. Put $$Q = \bigoplus_{j=1}^\infty p_{I_j} \in
  \cM(C(Z) \otimes \cK).$$ 
\begin{enumerate}
\item If $\big| \bigcup_{j \in F} I_j \big| \ge |F|$ for all finite
  subsets $F$ of $\N$, then $g \nprecsim Q$ and $Q$ is not properly infinite.
\item $g \precsim p_n \oplus p_n$ for every natural number $n$.
\item If infinitely many of the sets $I_1,I_2, \dots$ are singletons,
  then $Q \oplus Q$ is properly infinite and $Q \oplus Q \sim 1$ in
  $\cM(C(Z) \otimes \cK)$. 
\end{enumerate}
\end{proposition}

\begin{proof} (i). We show first that $g \nprecsim
  Q$ in $\cM(C(Z) \otimes \cK)$. Indeed, assume to the contrary that
  $g \precsim Q$. Then  
\begin{equation} \label{eq:f<Q}
g \precsim p_{I_1} \oplus p_{I_2} \oplus \cdots \oplus p_{I_k}
\end{equation}
in $C(Z) \otimes \cK$ for some $k$ by Lemma~\ref{lm:q<Q}. As noted
earlier, $C(Z) \otimes \cK$ is an inductive limit
$$\xymatrix{C(S^2) \otimes \cK \ar[rr]^-{\widehat{\rho}_1 \otimes
    \Id_\cK} && C((S^2)^2) \otimes \cK \ar[rr]^-{\widehat{\rho}_2
    \otimes \Id_\cK} && C((S^2)^3) \otimes \cK \ar[r] & \cdots \ar[r]
  & C(Z) \otimes \cK.}$$ 
Take $N$ such that all projections appearing in
\eqref{eq:f<Q} belong to the image of 
$$\widehat{\rho}_{\infty,n} \otimes \Id_\cK \colon C((S^2)^n) \otimes
\cK \to C(Z) 
\otimes \cK$$
whenever $n \ge N$. Use a standard inductive limit argument to see
that \eqref{eq:f<Q} holds 
relatively to $C((S^2)^n) \otimes \cK$ for some large enough $n \ge
N$. In the language  
of vector bundles over $(S^2)^n$, \eqref{eq:f<Q} and
Lemma~\ref{lm:Euler-p} imply that
\begin{equation} \label{eq:euler-obstruction}
\theta \oplus \eta \cong \zeta_{I_1} \oplus \zeta_{I_2} \oplus
\cdots \oplus \zeta_{I_k}
\end{equation}
for some vector bundle $\eta$ over $(S^2)^n$. Now,
\eqref{eq:euler-obstruction} and \eqref{eq:Euler1} imply that
$e(\zeta_{I_1} \oplus  
\cdots \oplus \zeta_{I_k})=0$, in contradiction with
Proposition~\ref{prop:marriage1} and the assumption on the sets $I_j$. 

The projection $p_{I_1}$ is a full element in $C(Z) \otimes \cK$ and
$p_{I_1} \le Q$. Hence $g$ belongs to the ideal generated by $Q$. It
now follows from Lemma~\ref{lm:prop} and from the fact that $ g
\nprecsim Q$ that $Q$ cannot be properly infinite.

(ii) follows from Lemma~\ref{lm:2zeta} and Lemma~\ref{lm:Euler-p}.

(iii). The unit $1$ of $\cM(C(Z) \otimes
\cK)$ can be written as a strictly convergent sum $1 =
\sum_{j=1}^\infty g_j$, where $g_j \sim g$ for all $j$. Let $\Gamma$
denote the infinite subset of $\N$ consisting of those $j$ for which
$I_j$ is a singleton. By Lemma~\ref{lm:rearranging} and (ii) we get
$$1 \; \sim \; \bigoplus_{j=1}^\infty g \; \precsim \; \bigoplus_{j \in
  \Gamma} (p_{I_j} \oplus p_{I_j}) \; \precsim  \; \bigoplus_{j=1}^\infty
(p_{I_j} \oplus p_{I_j}) \; \sim \; Q \oplus Q.$$
Lemma~\ref{lm:Qsim1} now tells us that $Q \oplus Q$ is properly
infinite and that $Q \oplus Q \sim 1$. 
\end{proof}

\section{A non-exact example} \label{sec:construction}
\noindent We construct here a simple, unital \Cs{} that
contains a finite and an infinite projection; thus proving one of our
main results: Theorem~\ref{thm:B} below. 

Let again $Z$ denote the infinite product space $\prod_{j=1}^\infty
S^2$. Set $A = C(Z) \otimes \cK = C(Z, \cK)$; recall from Section~4
that $\cM(A)$ 
denotes the multiplier algebra of $A$ and that it can be identified
with the set of bounded $^*$-strongly continuous functions 
$f \colon Z \to B(H)$.

Choose an injective function $\nu \colon \Z \times \N \to
\N$.  Choose points
$c_{j,i} \in S^2$ for all $j,i \in \N$ with $j \ge i$ such that
\begin{equation} \label{eq:c}
\overline{\{(c_{j,1},c_{j,2}, \dots, c_{j,n}) \mid j \ge n \}} = S^2
\times S^2 \times \cdots \times S^2
\end{equation}
for every natural number $n$. Set
\begin{equation} \label{eq:I_j}
I_j = \{\nu(j,1), \nu(j,2), \dots, \nu(j,j)\},
\end{equation}
for $j \in \N$.

Define \sh s $\varphi_j
\colon A \to A$ for all integers $j$ as follows. For $j \le 0$, set
\begin{equation} \label{eq:varphi1}
\varphi_j(f)(x) = f\big( x_{\nu(j,1)},x_{\nu(j,2)},x_{\nu(j,3)},
\dots\big), \qquad f \in A, \; x = (x_1,x_2, \dots ) \in Z.
\end{equation}
Let $p_n$ and $p_I$ be the projections in $A=C(Z,\cK)$ defined in
\eqref{eq:p_n} and \eqref{eq:p_I}. Choose an isomorphism $\tau \colon \cK
\otimes \cK \to \cK$.
For $f$ in $A$, $x=(x_1,x_2, \dots)$
in $Z$, and $j \ge 1$ define 
\begin{equation}
\label{eq:varphi2}
\varphi_j(f)(x) = \tau\big(f(c_{j,1}, \dots,
c_{j,j}, x_{\nu(j,j+1)}, x_{\nu(j,j+2)}, \dots ) \otimes
p_{I_j}(x)\big).
\end{equation}
Choose a sequence $\{S_j\}_{j=-\infty}^\infty$ of isometries in $\cM(A)$ such that
$\sum_{j=-\infty}^\infty S_jS_j^* = 1$ with the sum being strictly
convergent. Define a \sh{} $\psi \colon A \to \cM(A)$ by
\begin{equation}
\label{eq:varphi3}
\psi(f) = \sum_{j=-\infty}^\infty S_j\varphi_j(f)S_j^*, \qquad f \in A.
\end{equation}

\begin{lemma} \label{lm:F}
Let $\{e_n\}_{n=1}^\infty $ be an increasing approximate unit
for $A$. Then $\{\psi(e_n)\}_{n=1}^\infty$
converges strictly to a projection $F \in \cM(A)$, and $F$ is equivalent
to the identity 1 in $\cM(A)$.
\end{lemma}

\begin{proof} If $\psi(e_n)$ converges strictly to $F \in \cM(A)$ for
  \emph{some} approximate unit $\{e_n\}$ for $A$, then this conclusion
  will hold 
  for \emph{all} approximate units for $A$. We can therefore take
  $\{e_n\}_{n=1}^\infty $ to be the approximate unit given by $e_n(x)
  = \widehat{e}_n$, where $\{\widehat{e}_n\}_{n=1}^\infty$ is an
  increasing approximate unit for $\cK$. 

We show first that $\{\varphi_j(e_n)\}_{n=1}^\infty$ converges
strictly to a projection $F_j$ in $\cM(A)$ for each $j \in \Z$. Indeed,
since $\varphi_j(e_n) = e_n$ when $j \le 0$ it follows that
$\varphi_j(e_n) \to 1$ strictly; and so $F_j = 1$ when $j \le
0$. Consider next the case $j \ge 1$. Here we have $\varphi_j(e_n)(x) =
\tau(\widehat{e}_n \otimes p_{I_j}(x))$. Extend $\tau \colon \cK \otimes 
\cK \to \cK$ to a strongly continuous unital \sh{} $\overline{\tau}
\colon B(H \otimes H) \to B(H)$ and define $F_j$ in $\cM(A)$ by
$F_j(x) = \overline{\tau}(1 \otimes p_{I_j}(x))$ for $x \in Z$. Then $F_j$ is a
projection and $\{\varphi_j(e_n)\}_{n=1}^\infty$ converges strictly to
$F_j$.  

Now,
$$\psi(e_n) = \sum_{j=-\infty}^\infty S_j \varphi_j(e_n) S_j^* \;
\underset{n \to \infty}{\overset{\mathrm{strictly}}{\longrightarrow}} \; \sum_{j=-\infty}^\infty
S_jF_jS_j^* \eqdef F \in \cM(A),$$
As $1 = F_0 \sim S_0 F_0 S_0^* \le
F$ it follows from Lemma~\ref{lm:Qsim1} that $F \sim 1$ in 
$\cM(A)$.
\end{proof} 

\noindent Take an isometry $T$ in $\cM(A)$ with $TT^* = F$ (where $F$
is an in Lemma~\ref{lm:F}). Define
\begin{equation}
\label{eq:varphi4}
\varphi(f) = T^*\psi(f)T = \sum_{j=-\infty}^\infty
T^*S_j\varphi_j(f)S_j^*T, \qquad f \in A. 
\end{equation}
Then $\varphi \colon A \to \cM(A)$ is a \sh{} that maps an approximate
unit for $A$ into a sequence in $\cM(A)$ that converges strictly to
the identity in $\cM(A)$ (by Lemma~\ref{lm:F} and the choice of
$T$). It follows from \cite[Proposition 2.5]{Lan:modules} that 
$\varphi$ extends to a unital \sh{} $\overline{\varphi} \colon \cM(A)
\to \cM(A)$.

We collect below some properties of the \sh s $\varphi$ and
$\overline{\varphi}$. A subset of a \Cs{} $A$
is called \emph{full in $A$} if it is not contained in any proper
closed two-sided ideal in $A$.

\begin{proposition} \label{prop:varphi}
Let $p_1$ be the projection in $A$ defined in \eqref{eq:p_n} and
let $g$ be a constant 1-dimensional projection in $A = C(Z,\cK)$.
\begin{enumerate}
\item $\varphi(g) \sim 1$ in $\cM(A)$, and so $\varphi(f)$ is full in
  $\cM(A)$ for every full element $f$ in $A$.
\item If $f$ is a non-zero element in $\cM(A)$, then
  $\overline{\varphi}(f)$ does not belong to $A$, and
  $A\overline{\varphi}(f)$ is full in $A$. 
\item If $f$ is a non-zero element in $\cM(A)$, then
  $A\overline{\varphi}^{\, k}(f)$ is full
  in $A$ for every $k \in \N$.  
\item None of the projections 
  $\overline{\varphi}^{\,k}(p_1)$, $k \in \N$, are properly infinite
  in $\cM(A)$. 
\end{enumerate}
\end{proposition}

\noindent It follows immediately from (ii) that $\overline{\varphi}$
and $\varphi$ are injective, and that $\overline{\varphi}(\cM(A)) \cap
A = \{0\}$ and $\varphi(A) \cap A = \{0\}$.

The proof of Proposition~\ref{prop:varphi} is divided into a few
lemmas, the first of which (included for emphasis)
is standard and follows from the fact that any
closed two-sided ideal in $C(Z,\cK)$ is equal to
$C_0(U,\cK)$ for some open subset $U$ of $Z$. 

\begin{lemma} \label{lm:idealC(X)} 
Let $f$ be an element in $A=C(Z,\cK)$. Then $f$ is full in 
$A$ if and only if $f(x) \ne 0$ for all $x \in Z$.
\end{lemma}

\vspace{.3cm} \noindent {\bf{Proof of
    Proposition~\ref{prop:varphi}~(i):}} Observe first that
  $\varphi_j(g)=g$ for every $j \le 0$. Accordingly,
$$1 \; \sim \; \bigoplus_{j=-\infty}^0 g \; \sim \; \sum_{j=-\infty}^0
T^*S_j\varphi_j(g)S_j^*T \; \le \; \varphi(g) \quad \text{in} \; \cM(A).$$
This and Lemma~\ref{lm:Qsim1} imply that $\varphi(g) \sim
1$ and that $g$ is full in $\cM(A)$. If $f$
is any full element in $A$, then the closed two-sided ideal generated
by $\varphi(f)$ 
contains $\varphi(g)$ and therefore all of $\cM(A)$. This proves the
second claim in (i).
\hfill $\square$

\vspace{.3cm} \noindent {\bf{Proof of
    Proposition~\ref{prop:varphi}~(ii):}} Take a non-zero element $f$
  in $\cM(A)$. There is an element $a$ in $A$ such that $af \ne
  0$. The two claims in (ii) will clearly follow if we can show that
  $\overline{\varphi}(af) \notin A$ and that $A\overline{\varphi}(af)$
  is full in $A$, and we can therefore, upon replacing $f$ by $af$,
  assume that $f$ is a non-zero element in $A= C(Z,\cK)$.

There are $\delta > 0$, $r \in \N$, and non-empty open subsets $U_1,
  \dots, U_r$ of $S^2$ such that
\begin{equation} \label{eq:a}
x \in U_1 \times U_2 \times \cdots \times U_r \times S^2 \times
S^2 \times \cdots \implies \|f(x)\| \ge \delta.
\end{equation}
Use \eqref{eq:c} to find an infinite set $\Lambda$ of integers $j \ge
r$ such that  
\begin{equation} \label{eq:b}
(c_{j,1}, c_{j,2}, \dots, c_{j,r}) \in U_1 \times U_2 \times \cdots
\times U_r \quad \text{for all} \; \; j \in \Lambda.
\end{equation}
It follows from Lemma~\ref{lm:idealC(X)}, \eqref{eq:varphi2},
\eqref{eq:a}, and \eqref{eq:b} that $\|\varphi_j(f)\| \ge \delta$ and 
$\varphi_j(f)$ is full in $A$ for every $j$ in the infinite set
$\Lambda$. This entails 
that $\varphi(f) = \sum_{j=-\infty}^\infty T^*S_j \varphi_j(f)S_j^*T$
does not belong to $A$. (A strictly convergent sum
$\sum_{j=-\infty}^\infty a_j$ of pairwise orthogonal elements from $A$
belongs to $A$ 
if and only if $\lim_{j \to \pm \infty} \|a_j\| =0$.) The closed
two-sided ideal in $A$ generated by $A\varphi(f)$ contains the full
element $\varphi_j(f) = S_j^*T\varphi(f)T^*S_j$ and therefore all of
$A$ (for each---and hence
at least one---$j$ in $\Lambda$).  \hfill $\square$ 

\vspace{.3cm} \noindent {\bf{Proof of
    Proposition~\ref{prop:varphi}~(iii):}} This follows from
injectivity of $\overline{\varphi}$ and
Proposition~\ref{prop:varphi}~(ii). \hfill $\square$  

\vspace{.3cm} \noindent We proceed to prove
Proposition~\ref{prop:varphi}~(iv). 

\begin{lemma} \label{lm:J'} Let $J$ be a finite subset of
  $\N$ and let $j$ be an integer. Then $\varphi_j(p_J) \sim p_{\alpha_j(J)}$,
  where
\begin{equation} \label{eq:J'}
\alpha_j(J) = \begin{cases} \nu(j,J), &  j \le 0 \\
\nu(j,J \! \setminus \! \{1,2, \dots, j\}) \cup I_j,&  j \ge 1.
 \end{cases}
\end{equation}
We have in particular that $\nu(j,J) \subseteq \alpha_j(J)$ for all
finite subsets $J$ of $\N$ and for all $j \in \Z$.
\end{lemma}

\begin{proof} Write $J = \{t_1,t_2, \dots, t_k\}$, where $t_1 < t_2 < \cdots <
t_k$. We consider first the case where $j \le 0$. Then
\begin{eqnarray*}
\varphi_j(p_J)(x) & = & p_J(x_{\nu(j,1)}, x_{\nu(j,2)}, x_{\nu(j,3)},
\dots) \\ &= & p(x_{\nu(j,t_1)}) \otimes p(x_{\nu(j,t_2)}) \otimes \cdots
\otimes p(x_{\nu(j,t_k)}) \\ & = & p_{\nu(j,t_1)}(x) \otimes
p_{\nu(j,t_2)}(x) \otimes \cdots \otimes p_{\nu(j,t_k)}(x) \;
= \; p_{\nu(j,J)}(x),
\end{eqnarray*}
as desired. 

Suppose next that $j \ge 1$. Let $m$ be such that $t_{m-1} \le j <
t_m$ (with the convention $t_0=0$). Put
\begin{eqnarray*}
q(x) & = & p_J(c_{j,1}, \dots, c_{j,j}, x_{\nu(j,j+1)}, x_{\nu(j,j+2)},
  \dots) \\ & = & p(c_{j,t_1}) \otimes \cdots \otimes p(c_{j,t_{m-1}}) \otimes
  p(x_{\nu(j,t_m)}) \otimes \cdots \otimes p(x_{\nu(j,t_k)}) \\ 
& = & p(c_{j,t_1}) \otimes \cdots \otimes p(c_{j,t_{m-1}}) \otimes
p_{\nu(j,t_m)}(x) \otimes \cdots \otimes p_{\nu(j,t_k)}(x)  \\
& = & p(c_{j,t_1}) \otimes \cdots \otimes p(c_{j,t_{m-1}}) \otimes
p_{\nu(j,J \setminus \{1,2, \dots , j\})}(x).
\end{eqnarray*}
Thus $q \sim p_{\nu(j,J \setminus \{1,2, \dots , j\})}$. Now, 
$\varphi_j(p_J)(x) = \tau(q(x) \otimes p_{I_j}(x))$. Hence
$\varphi_j(p_J)$ is equivalent to the projection defined by 
$$x \mapsto \tau(p_{\nu(j,J \setminus \{1,2, \dots , j\})}(x) \otimes
p_{I_j}(x)),$$
and this projection is equivalent to $p_{\nu(j,J \setminus \{1,2,
  \dots , j\}) \cup I_j}$. 

The last claim follows from the definition of the sets $I_j$ in \eqref{eq:I_j}.
\end{proof}

\begin{lemma} \label{lm:marriage} Let $J_1,J_2, \dots$ be finite
  subsets of $\N$. Put $Q = \bigoplus_{i=1}^\infty p_{J_i} \in
  \cM(A)$. Then
$$\overline{\varphi}(Q)  \; \sim \; \bigoplus_{i=1}^\infty
\bigoplus_{j=-\infty}^\infty p_{\alpha_j(J_i)}.$$
Moreover, if $|\bigcup_{i \in F} J_i| \ge |F|$ for all finite subsets $F$
  of $\N$, then $|\bigcup_{(j,i) \in G} \alpha_j(J_i)| \ge |G|$ for
  all finite subsets $G$ of $\Z \times \N$. 
\end{lemma}

\begin{proof} By \eqref{eq:DS2}, $Q = \sum_{i=1}^\infty T_i p_{J_i}
  T_i^*$; and because $\overline{\varphi}$ is strictly continuous we get
$$\overline{\varphi}(Q) = \sum_{i=1}^\infty \overline{\varphi}(T_i)
\varphi(p_{J_i}) \overline{\varphi}(T_i)^* \sim 
\bigoplus_{i=1}^\infty \varphi(p_{J_i}) \sim
\bigoplus_{i=1}^\infty \bigoplus_{j=-\infty}^\infty \varphi_j(p_{J_i})
\sim \bigoplus_{i=1}^\infty
\bigoplus_{j=-\infty}^\infty p_{\alpha_j(J_i)},$$
where the first equivalence is proved below
\eqref{eq:DS1}--\eqref{eq:DS4}, and
the last equivalence follows from Lemma~\ref{lm:J'}.

By the Marriage Theorem we can find natural numbers
  $t_i \in J_i$ such that $\{t_i\}_{i \in \N}$ are mutually
  distinct. Set $s_{j,i} = \nu(j,t_i)$. Then $s_{j,i}$ belongs to
  $\alpha_j(J_i)$ by Lemma~\ref{lm:J'}, and $\{s_{j,i}\}_{(j,i) \in \Z
    \times \N}$ are mutually
  distinct because $\nu$ is injective and the $t_i$'s are mutually
  distinct. This proves the second claim of the lemma.
\end{proof}

\vspace{.3cm} \noindent {\bf{Proof of
    Proposition~\ref{prop:varphi}~(iv):}} Put $Q_0 = p_1$ and put
$Q_{n} = \overline{\varphi}^n(Q_0)$. We must show that none of the
projections $Q_n$, $n \ge 0$, are
  properly infinite. It is clear that $Q_0$ is finite, and
  hence not properly infinite.

Use Lemmas~\ref{lm:J'} and \ref{lm:marriage} to see that
$$Q_1 \; = \; \sum_{j=-\infty}^\infty T^*S_j\varphi_j(p_1)S_j^*T \;
\sim \; \bigoplus_{j=-\infty}^\infty 
\varphi_j(p_1) \; \sim \; \bigoplus_{j=-\infty}^0 p_{\nu(j,1)} \oplus
\bigoplus_{j=1}^\infty p_{I_j} \; = \; \bigoplus_{j=-\infty}^\infty p_{J_j},$$
where $J_j = \{\nu(j,1)\}$ for $j \le 0$ and $J_j = I_j$ for $j \ge
1$. It is easily seen that the sequence of sets $\{J_j\}_{j
  =-\infty}^\infty$ satisfies the condition $|\bigcup_{j \in F} J_j|
\ge |F|$ for all finite subsets $F$ of $\Z$. Hence $Q_1$ is not
properly infinite by  Proposition~\ref{prop:marriage2}~(i). 

The claim that $Q_n$ is not properly infinite for all $n$ follows by
induction using Lemma~\ref{lm:marriage} and
Proposition~\ref{prop:marriage2}~(i). \hfill $\square$

\begin{theorem} \label{thm:B}
Consider the inductive limit $B$ of the sequence:
$$\xymatrix@C-.5pc{\cM(C(Z) \otimes \cK) \ar[r]^-{\overline{\varphi}} & \cM(C(Z) \otimes \cK)
  \ar[r]^-{\overline{\varphi}} & \cM(C(Z) \otimes \cK)
  \ar[r]^-{\overline{\varphi}} & \cdots \ar[r] & B.}$$ 
Then $B$ has the following properties:
\begin{enumerate}
\item $B$ is unital and simple.
\item The unit of $B$ is infinite.
\item $B$ contains a non-zero finite projection.
\item $K_0(B) = 0$ and $K_1(B)=0$.
\end{enumerate}
\end{theorem}

\begin{proof} (i). $B$ is unital being the inductive limit of a
  sequence of unital \Cs s with unital connecting maps.

Write again $A$ for $C(Z) \otimes \cK$, and let
$\overline{\varphi}_{\infty,n} \colon \cM(A) \to B$ be the inductive
limit map from the $n$th copy of $\cM(A)$ into $B$. Let $L$ be a non-zero
closed two-sided ideal in $B$, and set 
$$L_n = \overline{\varphi}_{\infty,n}^{\;-1}(L) \vartriangleleft  \cM(A).$$
Then $L_n$ is non-zero for some $n$. Since $A$ is an essential ideal
in $\cM(A)$, 
also $A \cap L_n$ is non-zero. 

Take a non-zero element $e$ in $A \cap L_n$. Then
$\overline{\varphi}(e)$ belongs to $L_{n+1}$, hence $A
\overline{\varphi}(e) \subseteq L_{n+1}$, and so it follows
from Proposition~\ref{prop:varphi}~(ii) that $A \subseteq
L_{n+1}$. Take now a full element $f$ in $A \subseteq L_{n+1}$. Then
$\overline{\varphi}(f)$ belongs to $L_{n+2}$. It follows 
from Proposition~\ref{prop:varphi}~(i) that $\overline{\varphi}(f)$ is
full in $\cM(A)$ and therefore $L_{n+2}=\cM(A)$. Hence $L= B$,
and this shows that $B$ is simple.

(ii). This is clear because the unit of $\cM(A)$ is infinite.

(iii). As in the proof of Proposition~\ref{prop:varphi}~(iv), set $Q_0=p_1$ and
$Q_{n} = \overline{\varphi}^n(Q_0)$ for 
$n \ge 1$. Put $Q = \overline{\varphi}_{\infty,0}(Q_0) \in B$. It is
shown in Proposition~\ref{prop:varphi}~(ii) that $\overline{\varphi}$
is injective, which implies that $\overline{\varphi}_{\infty,0}$ is
injective, and hence $Q$ is non-zero. We show next that $Q$ is finite. 

Assume that $Q$ were infinite. Then $Q$ is properly infinite 
by Cuntz' result (see Proposition~\ref{prop:prop}) because $B$ is
simple. Applying 
Proposition~\ref{prop:inductive-limit} to the  sequence
$$
\xymatrix{
Q_0 \cM(A)Q_0 \ar[r]^{\lambda_0} & Q_1\cM(A)Q_1  \ar[r]^{\lambda_1} &
Q_2\cM(A)Q_2 \ar[r] &
\cdots  \ar[r] & QBQ, 
}$$
with the unital connecting maps $\lambda_j = \overline{\varphi}|_{Q_j
  \cM(A) Q_j}$, we obtain that $Q_n$ is properly infinite for all
sufficiently large 
$n$. But this contradicts Proposition~\ref{prop:varphi}~(iv).

(iv). This follows from the fact that the multiplier algebra of a
stable \Cs{} has trivial $K$-theory (see \cite[Proposition
12.2.1]{Bla:k-theory}).
\end{proof}

\vspace{.3cm} \noindent It follows
Proposition~\ref{prop:marriage2}~(ii) and
Proposition~\ref{prop:varphi}~(i) that the finite projection $Q$ in
$B$ (found in part (iii) above) satisfies
$$Q \oplus Q \; \sim \; \overline{\varphi}_{\infty,0}(Q_0 \oplus Q_0)
\; = \;
\overline{\varphi}_{\infty,0}(p_1 \oplus p_1) \; \succsim \;
\overline{\varphi}_{\infty,0}(g) =
\overline{\varphi}_{\infty,1}(\varphi(g)) \; \sim \; 1,$$
whence $Q \oplus Q \sim 1$ by Lemma~\ref{lm:Qsim1}. In other words,
the corner \Cs{} $QBQ$ is unital, finite, and simple, and $M_2(QBQ) \cong
B$ is infinite.

The \Cs{} $B$ from Theorem~\ref{thm:B} is not separable and not
exact. To see the latter, note that $B(H)$, the bounded operators on
a separable, infinite dimensional Hilbert space $H$, can be embedded
into $\cM(A) = \cM(C(Z) \otimes \cK)$ and hence into $B$. As $B(H)$ is
non-exact (see Wasserman \cite[2.5.4]{Was:exact}) it follows from Kirchberg's
result that exactness passes to sub-\Cs s (see
\cite[2.5.2]{Was:exact}) that $B$ is non-exact. We use the lemma
below from \cite{Bla:embed_simple} to construct a non-exact
\emph{separable} example. 

\begin{lemma}[Blackadar] \label{lm:Bla}
Let $B$ be a simple \Cs{} and let $X$ be a countable subset of $B$. It 
follows that $B$ has a separable, simple sub-\Cs{} $B_0$ that contains $X$.
\end{lemma}

\begin{corollary} \label{cor:B_0}
There exists a unital, \emph{separable}, non-exact, simple \Cs{} $B_0$
such that $B_0$ contains an infinite and a non-zero finite projection. 
\end{corollary}

\begin{proof} Let $B$ be as in Theorem~\ref{thm:B}. Let $s$ be a
  non-unitary isometry in $B$ and let $q$ be a non-zero finite
  projection in $B$. The universal \Cs{}, $C^*(\F_2)$,
  generated by two unitaries is separable and non-exact 
(see Wassermann \cite[Corollary~3.7]{Was:exact}). It admits an embedding into
$\cM(C(Z) \otimes \cK)$ and hence into $B$. Let $u,v \in B$ be the
images of the two 
(canonical) unitary generators in $C^*(\F_2)$. Use
Lemma~\ref{lm:Bla} to find a
separable, simple, and unital \Cs{} $B_0$ that contains $\{u,v,s,q\}$. 

Then $B_0$ is infinite because it
contains the non-unitary isometry $s$; and it contains the finite
projection $q$. Finally, $B_0$ is non-exact because it contains the
non-exact sub-\Cs{} $C^*(u,v) \cong C^*(\F_2)$. 
\end{proof}

\section{A nuclear example} \label{sec:addendum}

\noindent We show here that an elaboration of the construction in
Section~5 yields a \emph{nuclear} and separable example of a simple
\Cs{} with a finite and an infinite projection. 

The construction requires that we make a specific choice for 
the injective map $\nu \colon \Z \times \N \to \N$ from Section~5. 

Let $\{\Lambda_r\}_{r=0}^\infty$ be a partition of the set $\N$
such that $\Lambda_0 = \{1\}$ and such that $\Lambda_r$ is infinite for
each $r \ge 1$. For each $r \ge 1$ choose an injective map $\gamma_r \colon \Z
\times \Lambda_{r-1} \to \Lambda_r$ and define $\nu \colon \Z \times
\N \to \N$ by:
\begin{equation} \label{eq:nu}
\nu(j,t) = \gamma_r(j,t), \qquad r \in \N, \; t \in \Lambda_{r-1}, \;
j \in \Z. 
\end{equation}
Observe that 
\begin{equation} \label{eq:nu1}
t \in \Lambda_r \iff \nu(j,t) \in \Lambda_{r+1}, \qquad j \in \Z.
\end{equation}
To see that $\nu$ is injective assume that $\nu(j,t) =
\nu(i,s)$. Then $\nu(j,t) = \nu(i,s) \in \Lambda_r$ for some $r \ge
1$. Therefore both $s$ and $t$ belong to
$\Lambda_{r-1}$. Now, $\gamma_r(j,t) = \nu(j,t) =
\nu(i,s) = \gamma_r(i,s)$, which entails that $(j,t)=(i,s)$ by
injectivity of $\gamma_r$.

Let $\alpha_j$ be as defined in Lemma~\ref{lm:J'}
(wrt.\ the new choice of $\nu$). Let $\Gamma_0 \subseteq P(\N)$ be the family 
containing the one set $\{1\}$, and set 
$$\Gamma_{n+1} = \{\alpha_j(I) \mid I \in \Gamma_n, \; j \in \Z\}
\subseteq P(\N),$$
for $n \ge 0$. Set $\Gamma = \bigcup_{n=0}^\infty \Gamma_n$. Observe
that each $I \in \Gamma$ is a finite subset of $\N$. 

Put $Q_0 = p_1 \in A$ (cf.\ \eqref{eq:p_n}) and put $Q_{n} =
\overline{\varphi}^{ \,n}(Q_0) \in \cM(A)$ (where
$\overline{\varphi}$ is 
the endomorphism on $\cM(A)$ defined in
Section~5 above Proposition~\ref{prop:varphi}). It then follows by
induction from Lemma~\ref{lm:marriage} that  
\begin{equation} \label{eq:Q_n}
Q_n \sim \bigoplus_{I \in \Gamma_n} p_I, \qquad n \ge 0,
\end{equation}
when $p_I \in A$ is as defined in \eqref{eq:p_I}. 

\begin{lemma} \label{lm:1}
There is an injective function $t \colon \Gamma \to \N$
  such that $t(I) \in I$ for all $I \in \Gamma$. It follows in
  particular that
$$\big| \bigcup_{I \in F} I \big| \ge |F|$$
for all finite subsets $F$ of $\Gamma$.
\end{lemma}

\begin{proof} Define $t$ recursively on each $\Gamma_n$
as follows. For $n=0$ we set $t(\{1\}) = 1$. Assume that
  $t$ has been defined on $\Gamma_{n-1}$ for some $n \ge 1$. Then
  define $t$ on $\Gamma_n$ by $t(\alpha_j(I)) = \nu(j,t(I))$ for $I
  \in \Gamma_{n-1}$ and $j \in \Z$. It follows from Lemma~\ref{lm:J'} that
$$t(I) \in I \implies t(\alpha_j(I)) \in \alpha_j(I), \qquad I \in
\Gamma, \; j \in \Z.$$
It therefore follows by induction that $t(I) \in I$  for all $I \in \Gamma$. 

We show next that $t(I) \in \Lambda_n$ if $I \in \Gamma_n$. This is
clear for $n=0$. Let $n \ge 1$ and let $I \in \Gamma_n$ be given. Then
$I = \alpha_j(I')$ for some $I' \in \Gamma_{n-1}$ and some $j \in
\Z$. It follows that $t(I) = t(\alpha_j(I')) = \nu(j,t(I'))$. Hence
$t(I) \in \Lambda_n$ if $t(I') \in \Lambda_{n-1}$, cf.\
\eqref{eq:nu1}. Now the claim follows by induction on $n$.

We proceed to show that $t$ is injective. If $I,J \in
  \Gamma$ are such that $t(I) = t(J)$, then $t(I) = t(J) \in
  \Lambda_n$ for some $n$,  whence $I,J$ both belong to $\Gamma_n$. It
  therefore suffices to show that $t|_{\Gamma_n}$ is injective for
  each $n$. We prove this by induction on $n$. It is trivial that
  $t|_{\Gamma_0}$ is injective. Assume that $t|_{\Gamma_{n-1}}$ is
  injective for some $n \ge 1$. Let $I,J \in \Gamma_n$ be such that
  $t(I)=t(J)$. Then $I = \alpha_i(I')$
  and $J = \alpha_j(J')$ for some 
  $i,j \in \Z$ and some $I',J' \in \Gamma_{n-1}$, and
$$\nu(i,t(I')) =  t(\alpha_i(I')) = t(I) = t(J) 
= t(\alpha_j(J')) = \nu(j,t(J')).$$
Since $\nu$ is injective we deduce that $i=j$ and $t(I') =
t(J')$. By injectivity of $t|_{\Gamma_{n-1}}$ we obtain $I'=J'$, and
this proves that $I = J$. It has now been shown that $t|_{\Gamma_n}$
is injective, and the induction step is complete. 
\end{proof}

\noindent Let $g \in A = C(Z,\cK)$ be a constant 1-dimensional
projection, and let $Q_n$ be as defined above \eqref{eq:Q_n}.

\begin{lemma} \label{lm:2}
For each natural number $m$ we have
$$g \nprecsim Q_0 \oplus Q_1 \oplus \cdots \oplus Q_m \quad
\text{in} \; \cM(A).$$
\end{lemma}

\begin{proof} From \eqref{eq:Q_n} (and Lemma~\ref{lm:rearranging})
  we deduce that 
$$Q_0 \oplus Q_1 \oplus \cdots \oplus Q_n \; \sim \;
\bigoplus_{I \in \Gamma_0 \cup \cdots \cup \Gamma_n} p_I.$$
The claim of the lemma now follows from
Proposition~\ref{prop:marriage2}~(i) together with
Lemma~\ref{lm:1}. 
\end{proof}

\noindent As in Theorem~\ref{thm:B} consider the inductive limit
\begin{equation} \label{eq:B}
\xymatrix@C-.5pc{\cM(A) \ar[r]^-{\overline{\varphi}}
  & \cM(A) \ar[r]^-{\overline{\varphi}} & \cM(A)
  \ar[r]^-{\overline{\varphi}} & \cdots \ar[r] & B,}
\end{equation}
where $A= C(Z) \otimes \cK$. 
Let $\mu_{\infty,n} \colon \cM(A) \to B$ be the inductive limit map (from the
$n$th copy of $\cM(A)$) for $n \ge 0$. The endomorphism
$\overline{\varphi}$ on 
$\cM(A)$ extends to an automorphism $\alpha$ on $B$ that satisfies 
$\alpha(\mu_{\infty,n}(x)) = \mu_{\infty,n}(\overline{\varphi}(x))$
for $x \in \cM(A)$ and all $n \in \N$. (The inverse of $\alpha$ is
on the dense subset $\bigcup_{n=0}^\infty \mu_{\infty,n}(\cM(A))$ of
$B$ given by $\alpha^{-1}(\mu_{\infty,n}(x)) = \mu_{\infty,n+1}(x)$.)  

Put $A_0 = \mu_{\infty,0}(A) 
\subseteq B$, put $A_{n} = \alpha^n(A_0) \subseteq B$ for all $n \in
\Z$, and put
\begin{equation} \label{eq:D}
D_n = C^*(A_{-n}, A_{-n+1}, \dots, A_0, \dots A_{n-1}, A_n), \qquad D
= \overline{\bigcup_{n=1}^\infty D_n }. 
\end{equation}
It is shown in Lemma~\ref{lm:41} below that each $D_n$ is a type I
\Cs, and so the \Cs{} $D$ is an inductive limit of type I
algebras. In particular, $D$ is nuclear and belongs to the UCT
class $\cN$. Moreover, $D$ is 
$\alpha$-invariant (by construction). Observe that 
$A_{m-n} = \mu_{\infty,n}(\overline{\varphi}^{\, m}(A))$ for all
non-negative integers $m$ and $n$.

Put $Q = \mu_{\infty,0}(p_1) \; (= \mu_{\infty,n}(Q_n))$ in $D
\subseteq B$, and, as above, let $g \in A = C(Z,\cK)$ be a constant
1-dimensional projection. 

\begin{lemma} \label{lm:Q0} The following two relations hold in $D$ and in
  $B$:
\begin{enumerate} 
\item $\mu_{\infty,0}(g) \precsim Q \oplus Q$.
\item $\mu_{\infty,0}(g) \nprecsim \bigoplus_{j=-N}^N \alpha^j(Q)$ for
  all natural numbers $N$.
\end{enumerate}
\end{lemma}

\begin{proof} (i) follows immediately from
  Proposition~\ref{prop:marriage2}~(ii). 

(ii).  Assume, to reach a contradiction, that
$\mu_{\infty,0}(g)   \precsim  \sum_{j=-N}^N \alpha^j(Q)$ in $B$ 
(or in $D$) for some $N \in \N$. For $j \ge -N$ we have
$$\alpha^j(Q) = \alpha^j(\mu_{\infty,0}(Q_0)) =
\alpha^j(\mu_{\infty,N}(\overline{\varphi}^{\, N}(Q_0))) =
\mu_{\infty,N}(\overline{\varphi}^{\, {N+j}}(Q_0)).$$  
The relation $\mu_{\infty,0}(g)   \precsim  \sum_{j=-N}^N \alpha^j(Q)$
can therefore be rewritten as
$$\mu_{\infty,N}(\overline{\varphi}^{\, N}(g)) \; \precsim \;
\bigoplus_{j=0}^{2N} \mu_{\infty,N}(\overline{\varphi}^{\,
  {j}}(Q_0)) \quad \text{in} \; B.$$
By a standard property of inductive limits this entails that
$$\mu_{M,N}(\overline{\varphi}^{\, N}(g)) \; \precsim \;
\bigoplus_{j=0}^{2N} \mu_{M,N}(\overline{\varphi}^{\,
  {j}}(Q_0)) \quad \text{in} \; \cM(A),$$
for some $M \ge N$, or, equivalently, 
$$\overline{\varphi}^{\, M}(g) \; \precsim \;
\bigoplus_{j=0}^{2N} \overline{\varphi}^{\, {j+M-N}}(Q_0) \; = \;
\bigoplus_{j=M-N}^{N+M} \overline{\varphi}^{\,
  {j}}(Q_0) \; = \; \bigoplus_{j=M-N}^{N+M} Q_j \; \precsim \;
\bigoplus_{j=0}^{N+M} Q_j \quad \text{in} \; \cM(A).$$
Use now that $g \precsim \overline{\varphi}^{\, M}(g)$ (which holds because
$\varphi_j(g) = g$ for $j \le 0$, cf.\ \eqref{eq:varphi1}) to conclude that 
$g \; \precsim \; \bigoplus_{j=0}^{N+M} Q_j$ in $\cM(A)$, in
contradiction with Lemma~\ref{lm:2}.  
\end{proof}

\noindent
Let $C$ be an arbitrary unital \Cs{} and let $\gamma$ be an automorphism on
$C$.

Let $\cK$ denote the compact operators on $\ell^2(\Z)$ and
let $\{e_{i,j}\}_{i,j \in \Z}$ be a set of matrix units for
$\cK$. Define a unital injective \sh{} $\psi \colon C \to \cM(C \otimes \cK)$
and a unitary $U \in \cM(C \otimes \cK)$ by
$$\psi(c) = \sum_{n \in \Z} \gamma^n(c) \otimes e_{n,n}, \qquad U =
\sum_{n \in \Z} 1 \otimes e_{n,n+1}, \qquad c \in C,$$
(the sums converge strictly in $\cM(C \otimes \cK)$). It is easily
seen that
$$U\psi(c)U^* = \psi(\gamma(c)), \qquad c \in C,$$
so that $\psi$ extends to a representation $\widetilde{\psi} \colon C
\rtimes_\gamma \Z \to \cM(C \otimes \cK)$. The following standard
argument shows that the representation
$\widetilde{\psi}$ is faithful. 

Put $V_t = \sum_{n \in
  \Z} 1 \otimes t^{-n}e_{n,n} \in \cM(C \otimes \cK)$ for $t \in \T$,
and check that $V_t$ is a unitary element that satisfies $V_t\psi(c)V_t^*
= \psi(c)$ and $V_tUV_t^* = tU$ for all $t \in \T$. Let $E \colon C
\rtimes_\gamma \Z \to 
C$ be the canonical faithful conditional expectation, and define $F
\colon \image({\widetilde{\psi}}) \to \image({\widetilde{\psi}})$ by
$F(x) = \int_\T V_t x {V_t}^* \, dt$. Then $F(\widetilde{\psi}(x)) =
\psi(E(x))$ for all $x \in C \rtimes_\gamma \Z$. Now, if
$\widetilde{\psi}(x)=0$ for some positive element $x$ in $C
\rtimes_\gamma \Z$, then $\psi(E(x))=F(\widetilde{\psi}(x)) =0$,
whence $E(x)=0$ (by injectivity of $\psi$), and $x=0$ (because $E$ is
faithful).

\begin{lemma} \label{lm:Q1}
Let $C$ be a unital \Cs{} and let $\gamma$ be an automorphism on $C$. Suppose
that $p,q$ are projections in $C$ such that
\begin{enumerate}
\item $p \precsim \bigoplus_{j=1}^m q$ in $C$ for some natural number
  $m$, and
\item $p \nprecsim \bigoplus_{j=-N}^N \gamma^j(q)$ for all natural
  numbers $N$.
\end{enumerate}
Then $q$ is not properly infinite in $C \rtimes_\gamma
\Z$. 
\end{lemma}

\begin{proof} It suffices to show that $\psi(q)$ is not properly infinite in
  $\cM(C \otimes \cK)$. Assume, to reach a contradiction, that
  $\psi(q)$ is properly infinite in $\cM(C \otimes \cK)$. Then
  $\bigoplus_{j=1}^m \psi(q) \precsim \psi(q)$ by
  Proposition~\ref{prop:prop}. As $q 
  \otimes e_{0,0} \le \psi(q)$ we can use (i) to obtain 
$$p \otimes e_{0,0} \; \precsim \; \bigoplus_{j=1}^m q \otimes e_{0,0}
\; \le \; \bigoplus_{j=1}^m \psi(q) \; \precsim \; \psi(q) \; =
\; \sum_{j=-\infty}^\infty \gamma^j(q) \otimes e_{j,j}$$
in $\cM(C \otimes \cK)$. By Lemma~\ref{lm:q<Q} this entails that 
\begin{equation*}
p \otimes e_{0,0} \; \precsim \; \sum_{j=-N}^N \gamma^j(q)
\otimes e_{j,j} \quad \text{in} \; C \otimes \cK,
\end{equation*}
for some $N \in \N$, or, equivalently, that $p
\precsim  \bigoplus_{j=-N}^N \gamma^j(q)$ in $C$, in contradiction with
as\-sump\-tion (ii). 
\end{proof} 

\noindent Returning now to our specific \Cs{} $B$ from \eqref{eq:B},
Lemmas~\ref{lm:Q0} and \ref{lm:Q1} imply that:

\begin{lemma} \label{lm:Q}
The projection $Q = \mu_{\infty,0}(p_1)$ is not properly infinite in
$B \rtimes_\alpha \Z$.  
\end{lemma}

\begin{lemma} \label{lm:41}
The \Cs{} $D_n=C^*(A_{-n}, A_{-n+1}, \dots, A_0, \dots, A_n)$ is of type I
for each $n \in \N$.
\end{lemma}

\begin{proof} Note first that
\begin{equation} \label{eq:AnAm}
A_n A_m \subseteq A_{\min\{n,m\}}, \qquad n,m \in \Z.
\end{equation}
Indeed, we can assume without loss of generality that $n \le m$, and then
deduce
$$A_nA_m = \alpha^n(\mu_{\infty,0}(A\overline{\varphi}^{\, m-n}(A)))
\subseteq \alpha^n(\mu_{\infty,0}(A)) = A_n.$$
Since $A \cap \overline{\varphi}^{\, m-n}(A) = \{0\}$ when $n < m$, cf.\
Proposition~\ref{prop:varphi}~(ii), it follows also that 
\begin{equation} \label{eq:AnAm2}
A_n \cap A_m = \{0\}, \qquad n \ne m. 
\end{equation}
Use \eqref{eq:AnAm} to see that the \Cs{} $D_{m,n}$ generated by
$A_m,A_{m+1}, \dots, A_n$, for $m \le n$, is equal to
\begin{equation} \label{eq:D_nm}
D_{m,n} = A_m + A_{m+1} + \cdots + A_{n-1} + A_n.
\end{equation}
(To see that the right-hand side of \eqref{eq:D_nm} is norm closed,
use successively the fact that 
if $E$ is a \Cs{}, $I$ is a closed two-sided ideal in $E$, and $F$ is
a sub-\Cs{} of $E$, then $I+F$ is a sub-\Cs{} of $E$.)
It follows from \eqref{eq:AnAm}, \eqref{eq:AnAm2}, and \eqref{eq:D_nm}
that we have a decomposition series 
$$0 \; \vartriangleleft \; A_{-n}  \; \vartriangleleft \; D_{-n,-n+1}
\; \vartriangleleft \; D_{-n,-n+2} \; \vartriangleleft \; \cdots \;
\vartriangleleft \; D_{-n,n-1} \vartriangleleft \; D_{-n,n}=D_n$$
for $D_n$ and that each successive quotient is isomorphic to $A=C(Z)
\otimes \cK$. This proves that $D_n$ is a type I \Cs. 
\end{proof}

\begin{lemma} \label{lm:4}
The crossed product \Cs{} $D \rtimes_\alpha \Z$ contains an infinite
projection and a non-zero
projection which is not properly infinite. The \Cs{} $D$ has no
non-trivial $\alpha^n$-invariant closed two-sided ideal for any
non-zero integer $n$.  
\end{lemma}

\begin{proof} The projection $Q = \mu_{\infty,0}(p_1)$ belongs to $A_0
  = \mu_{\infty,0}(A) \subseteq D$, and it is non-zero
  because $\mu_{\infty,0}$ is injective (which again is because
  $\overline{\varphi}$ is injective). We have $D \subseteq B$ and
  hence
$$ Q \in D \rtimes_\alpha \Z \subseteq B \rtimes_\alpha \Z.$$
Since $Q$ is not properly infinite in $B \rtimes_\alpha \Z$ (by Lemma~\ref{lm:Q}) it
follows that $Q$ is not properly infinite in $D \rtimes_\alpha \Z$. 

Put $P = \mu_{\infty,0}(g) \in A_0 \subseteq D$, where $g$ is a
constant 1-dimensional projection in $A = C(Z,\cK)$. We have
$$g = \varphi_0(g) \sim S_0\varphi_0(g)S_0^* < \sum_{j=-\infty}^\infty
S_j\varphi_j(g)S_j^* = \overline{\varphi}(g),$$
cf.\ \eqref{eq:varphi1}. Hence $P =  \mu_{\infty,0}(g)$ is equivalent to
a proper subprojection of $\mu_{\infty,0}(\overline{\varphi}(g))$. As
$\mu_{\infty,0}(\overline{\varphi}(g))=\alpha(\mu_{\infty,0}(g)) \sim P$ in $D
\rtimes_\alpha \Z$ we conclude that
$P$ is an infinite projection in $D \rtimes_\alpha \Z$.

Suppose that $n$ is a non-zero integer (that we can take to be
positive) and that $I$ is a non-zero
closed two-sided $\alpha^n$-invariant 
ideal in $D$. Then $I \cap D_{kn}$ is non-zero for some natural number
$k$, cf.\ \eqref{eq:D}. As $I$ is
$\alpha^n$-invariant, $I \cap \alpha^{kn}(D_{kn})$ is non-zero,
and
$$\alpha^{kn}(D_{kn}) = C^*(A_0,A_1,\dots, A_{2kn}) =
\mu_{\infty,0}\big(C^*(A,\overline{\varphi}(A), \dots,
\overline{\varphi}^{\, 2kn}(A))\big).$$
Because $A_0 = \mu_{\infty,0}(A)$ is an essential ideal in
$\alpha^{kn}(D_{kn})$ it follows that 
$I \cap A_0$ is non-zero. Take a non-zero element $f$ in $I \cap A_0$,
and write $f = \mu_{\infty,0}(f_0)$ for some non-zero element $f_0$ in
$A$. Use Proposition~\ref{prop:varphi}~(iii) to conclude that
$$A_{-m}f = \mu_{\infty,m}\big(A
\overline{\varphi}^{\, m}(f_0) \big)$$
is full in $\mu_{\infty,m}(A)=A_{-m}$, and hence that $A_{-m}
\subseteq I$, for every natural number $m$. Since $I$ is
$\alpha^n$-invariant, $A_{-m+rn} = \alpha^{rn}(A_{-m}) \subseteq I$
for all $m \in \N$ and all $r \in \Z$. This shows that $A_m \subseteq
I$ for all $m$, which finally entails that $I=D$.  
\end{proof} 

\noindent We remind the reader of the notion of properly outer
automorphism introduced by Elliott in \cite{Ell:Auto}:

\begin{definition} An automorphism $\gamma$ on a \Cs{} $E$ is called
  \emph{properly outer} if for every non-zero $\gamma$-invariant
  closed two-sided ideal $I$ of $E$ and for every unitary $u$ in
  $\cM(I)$ one has $\|\gamma|_I - \Ad u \| = 2$ (the norm is the
  operator norm). 
\end{definition}

\noindent Olesen and Pedersen list in \cite[Theorem
6.6]{OlePed:C*-dynamicIII} eleven conditions on an automorphism
$\gamma$ that all are equivalent to $\gamma$ being properly outer. We
shall use the following sufficient (but not necessary) condition for
being properly outer: If
$E$ has no non-trivial $\gamma$-invariant ideals and if $\gamma(p)
\nsim p$ for some projection $p$ in $E$, then $\gamma$ is properly
outer. To see this, note first that $p \sim upu^* = (\Ad u)(p)$ for
every unitary $u$ in $\cM(E)$ (the equivalence holds relatively to
$E$). We therefore have $\gamma(p) \nsim (\Ad u)(p)$, whence
$\|\gamma(p)-(\Ad u)(p)\|=1$. This shows that $\|\gamma - \Ad u\| \ge
1$ for all unitaries $u$ in $\cM(E)$, whence $\gamma$ is properly outer
(by (ii) $\Leftrightarrow$ (iii) of \cite[Theorem
6.6]{OlePed:C*-dynamicIII}). 

(One can argue along another line by taking an approximate unit
$\{e_\lambda\}$ for $E$, such that $e_\lambda \ge p$ for all $\lambda$,
and set $x_\lambda = 2p-e_\lambda$. Then $x_\lambda$ is a contraction
in $E$ for all 
$\lambda$, and one can check that $\lim_{\lambda \to \infty}
\|\gamma(x_\lambda) - (\Ad u)(x_\lambda)\| = 2$, thus showing directly
that $\|\gamma - \Ad u\| = 2$ for all unitaries $u$ in $\cM(E)$
whenever $\gamma(p) \nsim p$ for some 
projection $p$ in $E$.) 

More generally, $\gamma$ is properly outer if for each
non-zero $\gamma$-invariant ideal $I$ of $E$ there is a projection $p$ in $I$
such that $\gamma(p) \nsim p$.

\begin{lemma} \label{lm:outer}
The automorphism $\alpha^n$ on $D$ is properly outer for every non-zero
integer $n$. 
\end{lemma}

\begin{proof} We know from Lemma~\ref{lm:4} that $D$ has no
  $\alpha^n$-invariant ideals (when $n \ne 0$), so the lemma will
  follow from the claim (verified below) 
  that $\alpha^n(Q) \nsim Q$ for all $n \ne 0$ (where $Q$ is as in
  Lemma~\ref{lm:Q0}). 

Assume, to reach a 
  contradiction, that $\alpha^n(Q) \sim Q$ for some non-zero integer
  $n$ (that we can take to be positive). Then, by
  Lemma~\ref{lm:Q0}~(i), $\mu_{\infty,0}(g)
  \precsim Q \oplus Q \sim Q \oplus \alpha^n(Q) \precsim \bigoplus_{j=0}^{n} 
  \alpha^j(Q)$ in $D$, in contradiction with
  Lemma~\ref{lm:Q0}~(ii).
\end{proof}

\noindent We now have all ingredients to prove our main result:

\begin{theorem} \label{thm:A}
There is a separable \Cs{} $D$ and an automorphism $\alpha$ on
$D$ such that:
\begin{enumerate}
\item $D$ is an inductive limit of type I \Cs s.
\item $D \rtimes_\alpha \Z$ is simple and contains an infinite and a
  non-zero finite projection.
\item $D \rtimes_\alpha \Z$ is nuclear and belongs to the UCT class $\cN$.
\end{enumerate}
\end{theorem} 

\begin{proof} Let $D$ be the \Cs{} and let $\alpha$ the automorphism on
  $D$ defined in (and above) \eqref{eq:D}. Since $D$
  is the union of an increasing sequence of sub-\Cs s $D_n$ (cf.\
  \eqref{eq:D}) and each $D_n$ is of type I (by Lemma~\ref{lm:41}), we
  conclude that $D$ is an inductive limit of type I \Cs s, and hence
  that the crossed product $D \rtimes_\alpha \Z$ is nuclear, separable,
  and belongs to the UCT class $\cN$.

Since $D$ has no non-trivial $\alpha$-invariant ideals (by 
Lemma~\ref{lm:4}) and $\alpha^n$ is properly outer for all $n \ne 0$
(by Lemma~\ref{lm:outer}), it follows from Olesen and Pedersen, \cite[Theorem
7.2]{OlePed:C*-dynamicIII}, (a result that extends results from
Elliott, \cite{Ell:Auto}, and Kishimoto, \cite{Kis:Auto}) that $D
\rtimes_\alpha \Z$ is simple. By simplicity of $D
\rtimes_\alpha \Z$, the (non-zero) projection $Q$, which in
Lemma~\ref{lm:4} is proved to be not properly infinite, must be finite in $D
\rtimes_\alpha \Z$, cf.\ Proposition~\ref{prop:prop}. The existence of
an infinite projection in $D \rtimes_\alpha \Z$ follows from
Lemma~\ref{lm:4}, and this completes the proof. 
\end{proof}

\section{Applications of the main results} \label{sec:main}

\noindent We begin by listing some corollaries to Theorems~\ref{thm:B}
and \ref{thm:A}. 

\begin{corollary} \label{cor:pi}
There is a nuclear, unital, separable, infinite, simple \Cs{} $A$ in
the UCT class $\cN$ such that $A$ is not purely infinite.
\end{corollary}

\begin{proof}
Take the \Cs{} $D \rtimes_\alpha \Z$ from Theorem~\ref{thm:A}, and take
a properly infinite projection $p$ and a non-zero finite projection $q$ in
that \Cs. Then $q \sim q_0 \le  p$ for some projection $q_0$ in 
$D \rtimes_\alpha \Z$ by Lemma~\ref{lm:prop}.
Hence $A = p(D \rtimes_\alpha \Z)p$
is infinite; and $A$ is not purely infinite because it contains the
non-zero finite projection $q_0$. 
\end{proof}

\begin{corollary} \label{cor:finite-traces}
There is a nuclear, unital, separable, finite, simple \Cs{} $A$ that is not
stably finite, and hence does not admit a tracial state (nor a non-zero
quasitrace). 
\end{corollary}

\begin{proof}
Take the \Cs{} $E=D \rtimes_\alpha \Z$ from Theorem~\ref{thm:A} and a non-zero
finite projection $q$ in $E$. Put $A = qEq$. Then $A$ is
finite, simple, and unital. Since $A \otimes \cK \cong E \otimes \cK$
we conclude that 
$A \otimes \cK$ (and hence $M_n(A)$ for some large enough $n$)
contains an infinite projection, so $A$ is not stably 
finite. 

Every simple, infinite 
\Cs{} is properly infinite, so $M_n(A)$ is properly infinite. No
properly infinite \Cs{} can admit a non-zero trace (or a quasitrace),
so $M_n(A)$, and hence $A$, do not admit a tracial state (nor a
non-zero quasitrace).  
\end{proof}

\noindent A \Cs{} $A$ is said to have the \emph{cancellation property} 
if the implication
\begin{equation} \label{eq:canc}
p \oplus r \sim q \oplus r \implies p \sim q
\end{equation}
holds for all projections $p,q,r$ in $A \otimes \cK$. It is known that
all \Cs s of stable rank one have the cancellation 
property and that no infinite \Cs{} has the cancellation
property. There is no example of a stably finite, simple \Cs{} which 
is known not to have the cancellation property (but 
Villadsen's \Cs s from \cite{Vil:sr=n} are
candidates). A \Cs{} $A$ is said to have the \emph{weak cancellation property} 
if \eqref{eq:canc} holds for those projections $p,q,r$ in $A \otimes
\cK$ where $p$ and
$q$ generate the same ideal of $A$.

\begin{corollary} \label{cor:perf}
There is a nuclear, unital, separable, simple \Cs{} $A$ that does not have the
weak cancellation property.
\end{corollary}

\begin{proof}
Take $A$ as in Corollary~\ref{cor:pi}, and take
a non-zero finite projection $q$ in $A$. Since $A$  
is properly infinite, we can find isometries
$s_1,s_2$ in $A$ with 
ortho\-gonal range projections; cf.\ Proposition~\ref{prop:prop}. Put $p
= s_1qs_1^*+ (1-s_1s_1^*)$. Then $p$ is 
infinite because $s_2s_2^* \le p$, and so $p \nsim q$ (because $q$ is
finite). On the other hand, $q$ and $p$ generate the same ideal of
$A$---namely $A$ itself---and 
\begin{eqnarray*}
p \oplus 1 & = & \big(s_1qs_1^*+ (1-s_1s_1^*) \big) \oplus 1 \; \sim
\; s_1qs_1^* \oplus  (1-s_1s_1^*) \oplus s_1s_1^*   \; \sim \; q \oplus 1.
\end{eqnarray*}
\end{proof}

\noindent It was  shown in \cite[Theorem 9.1]{KirRor:pi2} that the
following implications hold for any separable
\Cs{} $A$ and for any free filter $\omega$ on $\N$:
\begin{eqnarray*}
A \; \;  \text{is purely infinite} & \implies & A \; \; \text{is weakly purely infinite} \\
 & \iff &  A_\omega \;  \text{is traceless} \\ & \implies & A \; \; \text{is
   traceless,} \end{eqnarray*}
and the first three properties are equivalent for all \emph{simple} \Cs s $A$. 
(A \Cs{} is here said to be \emph{traceless} if no algebraic ideal in
$A$ admits a
non-zero quasitrace. See \cite{KirRor:pi2} for the definition of being 
weakly purely infinite.) It was not known in \cite{KirRor:pi2} if the
reverse of the third implication holds (for simple or for
non-simple \Cs s), but we can now answer this in the negative:

\begin{corollary} \label{cor:traces} Let $\omega$ be any free filter
  on $\N$. There is a nuclear, unital, separable,
  simple \Cs{} $A$ which is traceless, but where $\ell^\infty(A)$ and
  $A_\omega$ admit non-zero quasitraces defined on some (possibly non-dense)
  algebraic ideal. 
\end{corollary}

\begin{proof} Take $A$ as in Corollary~\ref{cor:finite-traces}. Then $A$ is
  algebraically simple and $A$ admits no (everywhere defined) non-zero
  quasitrace. Hence $A$ is traceless in the sense of
  \cite{KirRor:pi2}. Because $A$ is simple and not purely infinite, $A_\omega$
  cannot be traceless. Since $A_\omega$ is a quotient of
  $\ell^\infty(A)$, the latter \Cs{} cannot be traceless either.
\end{proof}

\noindent
Kirchberg has shown in \cite{Kir:fields} (see also \cite[Theorem
4.1.10]{Ror:encyc}) that every exact simple \Cs{} which
is \emph{tensorially non-prime} (i.e., is isomorphic to a tensor
product $D_1 \otimes D_2$, where $D_1$ and $D_2$ both are 
simple non-type~I \Cs s) is either stably finite or purely
infinite. Liming Ge has proved in \cite{Ge:prime}
that the II$_1$-factor ${\mathcal{L}}(\F_2)$ is
(tensorially) prime (in the von Neumann algebra sense), and it follows
easily from this result that the \Cs{}  
$C^*_{\mathrm{red}}(\F_2)$ is tensorially prime. We can now exhibit a simple,
\emph{nuclear} \Cs{} that is tensorially prime:

\begin{corollary} \label{cor:prime} The \Cs{} $D \rtimes_\alpha \Z$
  from Theorem~\ref{thm:A} is simple, separable, nuclear, and tensorially prime,
  and so is $p(D \rtimes_\alpha \Z)p$ for every non-zero projection
  $p$ in $D \rtimes_\alpha \Z$. 
\end{corollary}

\begin{proof} The \Cs{}  $D \rtimes_\alpha \Z$ is simple, separable,
  nuclear; cf.\ Theorem~\ref{thm:A}. It is not stably finite because
  it contains an infinite projection, and it is not purely infinite
  because it contains a non-zero finite projection. The (unital) \Cs{}
  $p(D \rtimes_\alpha \Z)p$ is stably isomorphic to $D \rtimes_\alpha
  \Z$ and is hence also simple, separable, nuclear, and neither stably
  finite nor purely infinite. It therefore follows from
  Kirchberg's theorem 
  (quoted above) that these \Cs s must be tensorially prime.
\end{proof}

\noindent Villadsen's \Cs s from \cite{Vil:perforation} and \cite{Vil:sr=n}
are, besides being simple and nuclear, probably also tensorially
prime (although to the knowledge of the author this has not yet been
proven). Jiang
and Su have in \cite{JiaSu:Z} found a non-type~I, unital, 
simple \Cs{} $\mathcal{Z}$ for which $A \cong A
\otimes \mathcal{Z}$ is known to hold for a large class of
well-behaved simple \Cs s $A$, such as for example the irrational
rotation \Cs s and more generally all \Cs s that are covered by a
classification theorem (cf.\ \cite{Ell:amenable?} or
\cite{Ror:encyc}). Such \Cs s $A$ are therefore not tensorially prime.

The real rank of the \Cs s found in Theorems~\ref{thm:B}
and \ref{thm:A} have not been determined, but we guess that they have real
rank $\ge 1$. That leaves open the following question:

\begin{question} \label{qu:rr0} Does there exist a (separable) unital,
  simple \Cs{} $A$ such that $A$ contains an infinite and a non-zero
  finite projection, and such that:
\begin{enumerate}
\item $A$ is of real rank zero?
\item $A$ is both nuclear and of real rank zero?
\end{enumerate}
\end{question}

\noindent It
appears to be difficult (if not impossible) to construct simple \Cs s
\emph{of real rank zero} that exhibit bad comparison properties; cf.\
Remark~\ref{rem:comp} below. 

George Elliott suggested the following:

\begin{question} \label{qu:inf} Does there exist a (separable),
  (nuclear), unital, simple \Cs{} $A$ such that all non-zero
  projections in $A$ are infinite but $A$ is not purely infinite?
\end{question}

\noindent If Question~\ref{qu:inf} has affirmative answer, and $A$ is
a unital, simple \Cs{} whose non-zero projections are infinite and $A$
is not purely infinite, then the real rank of $A$ cannot be zero. Indeed, a
simple \Cs{} is purely infinite if and only if it has real rank zero
and all its non-zero projections are infinite.

\begin{remark}[Comparison and dimension ranges] \label{rem:comp}
\noindent Suppose that $A$ is a unital, simple, infinite \Cs{} with a
non-zero finite projection $e$. By simplicity of $A$ there is a
natural number $k$ such that $1 \precsim e \oplus e \oplus \cdots
\oplus e$ (with $k$ copies of $e$).
Let $s_1,s_2, \dots$
be a sequence of isometries in $A$ with ortho\-gonal range
projections; cf.\ Proposition~\ref{prop:prop}. Letting $[p]$ denote
the Murray--von Neumann equivalence 
class of the projection $p$, we have 
$$n[1] = [s_1s_1^* + s_2s_2^* + \cdots + s_ns_n^*] \le [1] \le k[e]$$
for every natural number $n$. But $[1] \nleq [e]$ because $e$ is
finite and $1$ is infinite. 

This shows that if $A$ is a simple \Cs{} with a finite and an infinite
projection, then the semigroup ${\mathcal{D}}(A)$ of
Murray--von~Neumann equivalence classes of projections in $A \otimes
\cK$ is not weakly unperforated. 

(An ordered abelian semigroup 
$(S,+, \le)$ is said to be \emph{weakly unperforated} if 
$$\forall \, g,h \in S \; \forall \, n \in \N : ng < nh \implies g \le h.$$
The
order structure on ${\mathcal{D}}(A)$ is the algebraic order given by
$g \le h$ if and only if $h = g + f$ for some $f$ in ${\mathcal{D}}(A)$.)

Villadsen showed in \cite{Vil:perforation} that $K_0(A)$, and also
the semigroup ${\mathcal{D}}(A)$, of a \emph{simple, stably finite}
\Cs{} $A$ can 
fail to be weakly unperforated. The present article is a natural
continuation of Villadsen's work to the stably infinite case.

Let $(S,+)$ be an abelian semigroup with a zero-element $0$. An
element $g \in S$ 
is called \emph{infinite} if $g+x=g$ for some non-zero $x \in S$, and 
$g$ is called \emph{finite} otherwise. The sets of finite,
respectively, infinite elements in $S$ are denoted by
$S_{\mathrm{fin}}$ and $S_{\mathrm{inf}}$. One has $S = S_{\mathrm{fin}}
\amalg S_{\mathrm{inf}}$ and $S+ S_{\mathrm{inf}} \subseteq
S_{\mathrm{inf}}$, but the sum of two finite elements can be
infinite.

It is standard and easy to see that the finite and infinite elements
in the semigroup ${\mathcal{D}}(A)$ are given by
\begin{eqnarray*}
{\mathcal{D}}_{\mathrm{fin}}(A) & = & \{[f] : f \;  \text{is a finite
  projection in} \;  A\}, \\
{\mathcal{D}}_{\mathrm{inf}}(A) & = & \{[f] : f \;  \text{is an infinite
  projection in} \;  A\}.
\end{eqnarray*}

If $A$ is a simple \Cs{} that contains an infinite projection, then
the Grothendieck map $\gamma \colon {\mathcal{D}}(A) \to K_0(A)$
restricts to an isomorphism ${\mathcal{D}}_{\mathrm{inf}}(A) \to
K_0(A)$ as shown by Cuntz in \cite[Section 1]{Cuntz:KOn}. We can
therefore identify ${\mathcal{D}}_{\mathrm{inf}}(A)$ with $K_0(A)$,
in which case we can write
$${\mathcal{D}}(A) =  {\mathcal{D}}_{\mathrm{fin}}(A)\, \amalg \, K_0(A).$$
Note that $[0]$ belongs to ${\mathcal{D}}_{\mathrm{fin}}(A)$, and that
${\mathcal{D}}_{\mathrm{fin}}(A) = \{[0]\}$ if and only if all
non-zero projections in $A \otimes \cK$ are infinite. One can
therefore detect the existence of non-zero finite elements in $A \otimes
\cK$ from the semigroup ${\mathcal{D}}(A)$; and $K_0(A)$ contains all
information about ${\mathcal{D}}(A)$ if and only if all non-zero
projections in $A \otimes \cK$ are infinite.

In general, when $A$ is simple and contains both infinite and non-zero
finite projections, then  ${\mathcal{D}}_{\mathrm{fin}}(A)$ 
can be very complicated and large. One can show that
${\mathcal{D}}_{\mathrm{fin}}(B)$ is uncountable, when $B$ is
as in Theorem~\ref{thm:B}. We have no description of
${\mathcal{D}}(A)$, when $A = D \rtimes_\alpha \Z$ from
Theorem~\ref{thm:A}.

We remark finally, that if $A$ is simple and if $g$ is a non-zero element
in ${\mathcal{D}}_{\mathrm{fin}}(A)$, then $ng \in
{\mathcal{D}}_{\mathrm{inf}}(A)$ for some $n \in \N$. In other words,
${\mathcal{D}}_{\mathrm{inf}}(A)$ eventually absorbs all non-zero elements in
${\mathcal{D}}(A)$.
\end{remark}

\noindent 
The example found in Theorem~\ref{thm:A} provides a counterexample to
Elliott's classification conjecture (see for example
\cite{Ell:amenable?}) as it is formulated (by the author)
in \cite[Section 2.2]{Ror:encyc}. The conjecture asserts that 
\begin{equation} \label{eq:invariant}
\big(K_0(A), K_0(A)^+, [1_A]_0, K_1(A), T(A), r_A \colon T(A) \to
S(K_0(A))\big)
\end{equation}
is a complete invariant for unital, separable, nuclear, simple \Cs
s. If $A$ is stably infinite (i.e., if $A \otimes \cK$ contains an
infinite projection), then $K_0(A)^+ = K_0(A)$ and $T(A)=
\emptyset$. The Elliott invariant for unital, simple, stably infinite
\Cs s therefore degenerates to the triple $(K_0(A), [1_A]_0,
K_1(A))$. (We say that $(K_0(A), [1_A]_0,K_1(A))  \cong (G_0,g_0,G_1)$
if there are group isomorphisms $\alpha_0 \colon K_0(A) \to G_0$ and
$\alpha_1 \colon K_1(A) \to G_1$ such that $\alpha_0([1_A]_0) = g_0$.)

\begin{corollary} \label{cor:conjecture} 
There are two non-isomorphic nuclear, unital, separable, simple,
stably infinite \Cs s $A$ 
and $B$ (both in the UCT class $\cN$) such that  
$$(K_0(A), [1_A]_0,K_1(A)) \; \cong \; (K_0(B), [1_B]_0,K_1(B)).$$
\end{corollary}

\begin{proof} Take the \Cs{} $A$ from Corollary~\ref{cor:pi}.  
It follows from \cite[Theorem 3.6]{Ror:infsimple} that there is a 
nuclear, unital, separable, simple, purely infinite \Cs{} $B$ in the
UCT class $\cN$ such that 
$$(K_0(A), [1_A]_0,K_1(A)) \; \cong \; (K_0(B), [1_B]_0,K_1(B)).$$
Since $B$ is purely infinite and $A$ is not purely infinite, we have
$A \ncong B$.  
\end{proof}

\noindent One can amend the Elliott invariant by replacing the triple
$(K_0(A), K_0(A)^+, [1_A]_0)$ (for a unital \Cs{} $A$) with the pair
$(\mathcal{D}(A), [1_A])$, cf.\ Remark~\ref{rem:comp} above, where
$\mathcal{D}(A)$ carries the structure of a semigroup. In the unital, stably
infinite case, the amended invariant will then become
$(\mathcal{D}(A), [1_A], K_1(A))$.
(Since $K_0(A)$ is the Grothendieck group of $\mathcal{D}(A)$, and
$K_0(A)^+$, respectively, $[1_A]_0$, are the images of
$\mathcal{D}(A)$, respectively, $[1_A]$, under the Grothendieck map
$\gamma \colon \mathcal{D}(A) \to K_0(A)$, 
one can recover $(K_0(A), K_0(A)^+, [1_A]_0)$ from $(\mathcal{D}(A),
[1_A])$.)

The invariant $(\mathcal{D}(A),
[1_A])$ can detect if $A$ has a non-zero finite projection,
cf.\ Remark~\ref{rem:comp}; and the triples $(\mathcal{D}(A), [1_A],
K_1(A))$ and $(\mathcal{D}(B), [1_B], K_1(B))$ are therefore
non-isomorphic, when $A$ and $B$ are as in
Corollary~\ref{cor:conjecture}. We have no 
example to show that $(\mathcal{D}(A),
[1_A], K_1(A))$ is not a complete invariant for nuclear, unital, simple,
separable, stably infinite \Cs s. On the other hand, there is no
evidence to suggests that
$(\mathcal{D}(A), [1_A], K_1(A))$ indeed is a complete invariant for
this class of \Cs s. 

The Elliott conjecture can also be amended by restricting the class of
\Cs s that are to be classified. One possibility is to consider only
those unital, separable, nuclear, simple \Cs s $A$ for which $A \cong
A \otimes \mathcal{Z}$ where $\mathcal{Z}$ is the Jiang--Su algebra
(see the comment 
below Corollary~\ref{cor:prime}). It seems plausible that the Elliott
invariant \eqref{eq:invariant} actually is a complete invariant for
this class of \Cs s; and one could hope that the condition $A \cong A
\otimes \mathcal{Z}$ has an alternative intrinsic equivalent
formulation, for example in terms of the existence of sufficiently many
central sequences.

\begin{remark}[A non-simple example] \label{rm:ext} Examples of 
  non-simple unital \Cs s $A$, such that $A$ is finite and $M_2(A)$ 
  is infinite, have been known for a long time. Such examples were
  independently discovered by Clarke in \cite{Cla:finite} and by
  Blackadar (see Blackadar \cite[Exercise 6.10.1]{Bla:k-theory}): One
  such example is obtained by taking a unital extension
$$\xymatrix{ 0 \ar[r] & \cK \ar[r] & A \ar[r] & C(S^3) \ar[r] & 0 }$$
with non-zero index map $\delta \colon K_1(C(S^3)) \to K_0(\cK)$. Then 
$A$ is finite and $M_2(A)$ is infinite. 

The proof uses that any
isometry or co-isometry $s$ in $A$ (or in a matrix algebra over $A$)
is mapped to a unitary  
element $u$ in (a matrix algebra over) $C(S^3)$; and every unitary $u$
in
$M_n(C(S^3))$ lifts to an isometry or a co-isometry
$s$ in $M_n(A)$. Moreover, the isometry or co-isometry $s$ is
non-unitary if and only if the unitary element $u$ has non-zero
index. The unitary group of $C(S^3)$ is connected, so all unitaries
here have zero index. Hence $A$ contains no non-unitary isometry, so
$A$ is finite. By construction of the extension, the
generator of $K_1(C(S^3))$, which is a unitary element in $M_2(C(S^3))$, has 
non-zero index, and so it lifts to a non-unitary isometry or co-isometry
in $M_2(A)$, whence $M_2(A)$ is infinite. 

The \Cs{} $M_2(A)$ is not properly infinite since the quotient,
$M_2(A)/M_2(\cK) \cong M_2(C(S^3))$, is finite.
\end{remark}

\noindent An example of a unital, finite, (non-simple) \Cs{} $A$ such
that $M_2(A)$ is properly infinite was found in \cite{Ror:sums}.

\begin{remark}[Inductive limits] \label{rem:inductive}
Suppose that
$$\xymatrix{B_1 \ar[r] & B_2 \ar[r] & B_3 \ar[r] & \cdots \ar[r] &
  B}$$
is an inductive limit with unital connecting maps, and that $B$ is a
simple \Cs{} such that $B$ is finite and $M_2(B)$ is infinite. Then $M_2(B)$
is properly infinite, and it follows from Proposition~\ref{prop:inductive-limit}
that $B_n$ is finite and $M_2(B_n)$ is properly infinite for all
sufficiently large $n$. It is therefore not possible to construct an
example of a simple \Cs{}, which is finite, but not stably finite, by
taking an inductive limit of \Cs s arising as in the example described
in Remark~\ref{rm:ext}.  
\end{remark}

\begin{remark}[Free products] \label{rm:free}
Let $B$ be a simple, unital \Cs{} such that $B$ is finite and $M_2(B)$ 
is infinite. Then we have unital \sh s
$$\varphi_1 \colon M_2(\C) \to M_2(B), \qquad \varphi_2 \colon
\cO_\infty \to M_2(B),$$
such that $\varphi_1(e)$ is a finite projection in $M_2(B)$ whenever $e$
is a one-dimensional projection in $M_2(\C)$. 

The existence of $B$ (already obtained in the non-simple case 
in \cite{Ror:sums}) shows that the image of $e$ in the universal
unital free product \Cs{} $M_2(\C) * \cO_\infty$ is not properly
infinite. 

It is tempting to turn this around and seek a simple \Cs{} $A$ with a
finite and an infinite projection by defining $A$ to be a suitable
free product of $M_2(\C)$ and $\cO_\infty$. However, the universal
unital free product $M_2(\C) * \cO_\infty$ is not simple. The reduced
free product \Cs{}
$$(A, \rho) = (M_2(\C), \rho_1) * (\cO_\infty, \rho_2),$$
with respect to faithful states $\rho_1$ and $\rho_2$, is simple (at least for
many choices of the states $\rho_1$ and $\rho_2$, see for example
\cite{Avi:simple}) and properly infinite, but no non-zero projection $e$ in
$M_2(\C)$ is finite in $A$. The Cuntz algebra $\cO_\infty$ contains a
sequence of non-zero mutually orthogonal projections, and it 
therefore contains a projection $f$ with $\rho_2(f) < \rho_1(e)$. Now,
$e$ and $f$ are free 
with respect to the state $\rho$ and 
$\rho(f) < \rho (e)$. This implies that $f \precsim e$ (see
\cite{ABH:choi}), and therefore $e$ must be infinite.

It is shown in \cite{DykRor:projections} that reduced free
product \Cs s often have weakly unperforated $K_0$-groups, which is
another reason why this class of \Cs s is unlikely to provide an
example of a simple \Cs{} with finite and infinite projections; cf.\
Remark~\ref{rem:comp}.  
\end{remark}

\noindent We conclude this article by remarking that ring theorists for a 
long time have known about finite simple \emph{rings} that are not stably finite:

\begin{remark}[An example from ring theory] \label{rm:ring} A unital
  ring $R$ is called \emph{weakly finite} if $xy=1$ implies $yx=1$ for all $x,y$ in
  $R$, and $R$ is called \emph{weakly $n$-finite} if $M_n(R)$ is weakly
  finite. (A finite ring is a ring with finitely many elements!)  
A (unital) non-weakly finite simple ring $R$ is properly infinite in the sense
  that there are idempotents $e,f$ in $R$ such that $1 \sim e \sim f$
  and $ef=fe=0$. (Equivalence of idempotents is given by $e \sim f$ if 
  and only if $e = xy$ and $f = yx$ for some $x,y$ in $R$.) 

An example of a unital, simple ring which is weakly finite but not
weakly 2-finite was constructed by P.\ M.\ Cohn as follows: 

Take natural numbers $2 \le m < n$ and consider the universal ring
$V_{m,n}$ generated by $2mn$ elements $\{x_{ij}\}$
and $\{y_{ji}\}$, $i = 1, \dots, m$ and $j=1, \dots, n$, satisfying the
relations $XY= I_m$ and $YX=I_n$, where $X = (x_{ij}) \in M_{m,n}(R)$, 
$Y = (y_{ij}) \in M_{n,m}(R)$, and $I_m$ and $I_n$ are the units of the 
matrix rings $M_m(R)$ and $M_n(R)$. The rings $M_m(V_{m,n})$ and $M_n(V_{m,n})$
are isomorphic and $M_n(V_{m,n})$ is not weakly finite. Therefore
$M_m(V_{m,n})$ is not weakly finite. In other words, $V_{m,n}$ is
not weakly $m$-finite.

It is shown by Cohn in \cite[Theorem 2.11.1]{Cohn:free} (see also the
remarks at the end of Section~2.11 of that book) that $V_{m,n}$ is a
so-called \emph{$(m-1)$-fir}, and hence a $1$-fir; and a ring is a
$1$-fir if and only if it is an integral domain (i.e., if it has no
non-zero zero-divisors). Cohn proved in \cite{Cohn:simple} that every integral
domain embeds into a simple integral domain. In particular, $V_{m,n}$
is a subring of a simple integral domain $R_{m,n}$ whenever $2 \le m < 
n$. Now, $R_{m,n}$ is weakly finite (an integral domain has no
idempotents other than $0$ and $1$ and must hence be weakly finite),
and $R_{m,n}$ is not weakly $m$-finite (because it contains
$V_{m,n}$). 

This example cannot in any obvious way be carried over to \Cs s, first of all 
because no \Cs{} other that $\C$ is an integral domain.
\end{remark}

{\small{
\providecommand{\bysame}{\leavevmode\hbox to3em{\hrulefill}\thinspace}

}}
\vspace{.4cm}

\noindent{\sc Department of Mathematics, University of Southern
  Denmark, Odense, Campusvej 55, 5230 Odense M, Denmark}

\vspace{.2cm}

\noindent{\sl E-mail address:} {\tt mikael@imada.sdu.dk}\\
\noindent{\sl Internet home page:}
{\tt www.imada.sdu.dk/$\,\widetilde{\;}$mikael/}
\end{document}